\documentclass[11pt]{amsart}
\usepackage{amsmath,amscd,amssymb,amsfonts,latexsym,graphics,color}
\usepackage[all]{xy}
\usepackage{ulem}
\usepackage[latin1]{inputenc}

\usepackage[T1]{fontenc}
\usepackage{eucal} 
\usepackage{amsbsy}
\usepackage{amsthm,alltt,dsfont}

\usepackage{geometry}
\allowdisplaybreaks

\newcommand{\bC}{{\mathbb C}}

\def\bF{\mathbb F}

\newcommand{\bK}{{\mathbb K}}
\newcommand{\bL}{{\mathbb L}}
\newcommand{\bM}{{\mathbb M}}
\newcommand{\bN}{{\mathbb N}}

\newcommand{\bQ}{{\mathbb Q}}
\newcommand{\bZ}{{\mathbb Z}}

\newcommand{\cM}{{\mathcal M}}

\newcommand{\cO}{{\mathcal O}}

\newcommand{\cV}{{\mathcal V}}

\newcommand{\cW}{{\mathcal W}}
\newcommand{\cY}{{\mathcal Y}}

\newcommand{\mh}{\mbox{MHM}}

\newcommand{\ms}{{\it mHs}}

\newcommand{\lra}{\longrightarrow}

\newcommand{\Aut}{{\it Aut}}

\theoremstyle{plain}
\newtheorem{thm}{Theorem}[section]
\newtheorem{cor}[thm]{Corollary}

\newtheorem{lem}[thm]{Lemma}
\newtheorem{prop}[thm]{Proposition}

\newtheorem*{ack}{Acknowledgements}
\theoremstyle{definition}
\newtheorem{df}[thm]{Definition}
\newtheorem{rem}[thm]{Remark}
\newtheorem{example}[thm]{Example}

\def\be{\begin{equation}}
\def\ee{\end{equation}}

\def\bt{\begin{thm}}
\def\et{\end{thm}}

\def\bc{\begin{cor}}
\def\ec{\end{cor}}

\def\br{\begin{rem}}
\def\er{\end{rem}}

\def\bp{\begin{prop}}
\def\ep{\end{prop}}

\def\bl{\begin{lem}}
\def\el{\end{lem}}

\def\bn{\begin{enumerate}}
\def\en{\end{enumerate}}

\def\bex{\begin{example}}
\def\eex{\end{example}}

\def\bd{\begin{df}}
\def\ed{\end{df}}

\title{Plethysm and cohomology representations of external and symmetric products}
\author{Lauren\c{t}iu Maxim}
\address{L. Maxim: Department of Mathematics, University of Wisconsin-Madison,  480 Lincoln Drive, Madison, WI 53706, 
USA.}
\email {maxim@math.wisc.edu}
\author[J. Sch\"urmann ]{J\"org Sch\"urmann}
\address{J.  Sch\"urmann : Mathematische Institut,
          Universit\"at M\"unster,
          Einsteinstr. 62, 48149 M\"unster,
          Germany.}
\email {jschuerm@math.uni-muenster.de}

\date{\today}


\begin{document}

\begin{abstract}  
We prove refined generating series formulae for characters of (virtual) cohomology representations of external products of suitable coefficients, e.g., (complexes of) constructible or coherent sheaves, or (complexes of) mixed Hodge modules on spaces such as (possibly singular) complex quasi-projective varieties. These formulae generalize our previous results for symmetric and alternating powers of such coefficients, and apply also to other Schur functors.  The  proofs of these results are reduced via an equivariant K\"{u}nneth formula to a more general generating series identity for abstract characters of tensor powers $\cV^{\otimes n}$  of an element $\cV$ in a suitable symmetric monoidal category $A$. 
This abstract approach applies directly also in the equivariant context for spaces with additional symmetries (e.g., finite group actions, finite order automorphisms, resp., endomorphisms), as well as for introducing an abstract plethysm calculus for symmetric sequences of objects in $A$.
\end{abstract}

\maketitle

\tableofcontents

\section{Introduction}

In this paper we consider ``nice'' spaces such as: compact topological spaces with finite dimensional  cohomology (e.g., compact complex analytic spaces), semi-algebraic sets, complex quasi-projective varieties, or varieties over any base field of characteristic zero.
In fact, in Section \ref{abssec} we explain our results from an abstract axiomatic point of view based on the equivariant K\"{u}nneth formula, 
which also covers cases like cohomology of rational Chow motives (\cite{He, RT}), or Zeta functions of constructible sheaves for the Frobenius endomorphism of varieties over finite fields
(as in \cite{SGA5}[Thm. on p.464] and \cite{FK}[Thm.4.4 on p.174]). 

\subsection{Generating series formulae}\label{gs}
In this paper, we obtain refined generating series formulae for characters of (virtual) cohomology representations of external products of  ``nice'' spaces $X$ as above. In fact, we only need a suitable cohomology theory $H^*$ (e.g., closed or compactly supported singular cohomology of topological spaces, or a Weil cohomology for smooth projective varieties over finite fields), so that $H^*(X)$ is a finite dimensional $\mathbb{K}$-vector space (with $\mathbb{K}$ a field of characteristic zero), and which satisfies a $\Sigma_n$-equivariant K\"unneth formula: 
\be\label{Kut} H^*(X^n) \simeq H^*(X)^{\otimes n},\ee
with $\Sigma_n$ the symmetric group on $n$-elements.

Our results and techniques also apply to  suitable coefficients $\cM$ on such a space $X$,  provided that $H^*_{({c})}(X,\cM)$ is a  finite dimensional $\mathbb{K}$-vector space, and a corresponding $\Sigma_n$-equivariant K\"unneth formula holds:
\be\label{Kue} H^*_{({c})}(X^n,\cM^{\boxtimes n}) \simeq H^*_{({c})}(X,\cM)^{\otimes{n}},\ee
where $\cM^{\boxtimes n}$ is the $n$-th self-external product of $\cM$ with its induced $\Sigma_n$-action. (Here $H^*_{({c})}$  denotes the (compactly supported) cohomology or hypercohomology.)
Examples of such coefficients include: 
\begin{enumerate}
\item[(a)] if $X$ is a locally compact topological space, we work with $\cM \in D_f^b(X;\mathbb{K})$, a bounded sheaf complex of $\mathbb{K}$-vector spaces, such that   $H^*_{c}(X,\cM)$ is  {\it finite} dimensional (in this case, the K\"unneth formula is available only for the compactly supported cohomology).
\item[(b)] if $X$ is a complex algebraic or real semi-algebraic set, respectively, a compact complex analytic, subanalytic, or Whitney stratified set, we work with $\cM \in D^b_c(X;\mathbb{K})$, a bounded complex of sheaves of $\mathbb{K}$-vector spaces, which is constructible (in the corresponding sense). 
Here constructibility also includes the assumption that all stalks are finite dimensional.
\item[({c})] if $X$ is a compact complex analytic set or a projective algebraic variety over a field $\mathbb{K}$ of characteristic zero, we consider $\cM \in D^b_{coh}(X)$, a bounded complex of $\cO_X$-modules with coherent cohomology. 
\item[(d)] if $X$ is a complex quasi-projective variety, we can work with $\cM \in D^b\mh(X)$, a bounded complex of algebraic mixed Hodge modules on $X$.
\end{enumerate}
For further discussions of the K\"unneth formula in the cases (a)-({c}) above, see \cite{MS}[Sect.4.4], while for the case (d) of mixed Hodge module coefficients see \cite{MSS} .
 
In all cases (a)-(d) above, the (compactly supported) cohomology $H^*_{({c})}(X^n,\cM^{\boxtimes n})$ is a $\Sigma_n$-representation. Let $Rep_{\mathbb{K}}(\Sigma_n)$ be the Grothendieck group of (finite dimensional) $\mathbb{K}$-representations of $\Sigma_n$. By associating to a representation its character, we get a group monomorphism (with finite cokernel):
$$tr_{\Sigma_n}:Rep_{\mathbb{K}}(\Sigma_n) \hookrightarrow C(\Sigma_n),$$ 
with $C(\Sigma_n)$ the free abelian group of $\bZ$-valued class functions on $\Sigma_n$. Recall that characters of a symmetric group are integer valued. Moreover, $Rep_{\mathbb{K}}(\Sigma_n)$ is freely generated by the irreducible $\bK$-representations, all of which come from rational representations.
 
Consider the generating {\it Poincar\'e polynomial for the characters} of the above $\Sigma_n$-representations, namely: 
$$tr_{\Sigma_n}(H^*_{({c})}(X^n,\cM^{\boxtimes n})):=\sum_{k}tr_{\Sigma_n}(H_{({c})}^{k}(X^n,\cM^{\boxtimes n})) \cdot (-z)^k \in C(\Sigma_n) \otimes_{\bZ} \bZ[z^{\pm 1}],$$
and similarly for the characters of the finite dimensional cohomology $H^*(X^n)$ in the context of (\ref{Kut}). 
Aditionally, in the case when $\cM \in D^b\mh(X)$, the cohomology groups $H^*_{({c})}(X^n,\cM^{\boxtimes n})$ carry mixed Hodge structures, and the associated graded complex vector spaces
$$H_{({c})}^{p,q,k}(X^n,\cM^{\boxtimes n}):=Gr^p_F Gr^W_{p+q}H_{({c})}^{k}(X^n,\cM^{\boxtimes n})$$ of the Hodge and resp. weight filtrations are also $\Sigma_n$-representations. So in this case we can also consider the following more refined generating {\it mixed Hodge polynomial for the characters} of the $\Sigma_n$-representations of these associated graded vector spaces, namely:
$$tr_{\Sigma_n}(H^*_{({c})}(X^n,\cM^{\boxtimes n})):=\sum_{p,q,k}tr_{\Sigma_n}(H_{({c})}^{p,q,k}(X^n,\cM^{\boxtimes n})) \cdot y^px^q(-z)^k \in C(\Sigma_n) \otimes_{\bZ} \bZ[y^{\pm 1},x^{\pm 1},z^{\pm 1}].$$
While we use the same notation for the two types of generating polynomials (Poincar\'e and, resp., mixed Hodge), the reader should be able to distinguish their respective meaning from the context.
Note that by forgetting the grading with respect to the mixed Hodge structure (i.e., by letting $y=x=1$), the mixed Hodge polynomial (defined for mixed Hodge module coefficients) specializes to the Poincar\'e polynomial for the underlying constructible sheaf complex.

To simplify the notations and statements even further, we let $\mathbb{L}$ denote any of the two Laurent polynomial rings $\bZ[z^{\pm 1}]$ and, respectively, $\bZ[y^{\pm 1},x^{\pm 1},z^{\pm 1}]$. Once again, its meaning in the results below should be clear from the context.

\medskip

In this paper, for a pair $(X,\cM)$ of a space $X$ with coefficient $\cM$ as in the cases (a)-(d) above, we aim to calculate 
the generating series:
$${\sum_{n \geq 0} tr_{\Sigma_n}(H^*_{({c})}(X^n,\cM^{\boxtimes n})) \cdot t^n} \in \bigoplus_n C(\Sigma_n) \otimes_{\bZ} \mathbb{L}[[t]]$$ in terms of the corresponding {\it Poincar\'e polynomial} $$P_{({c})}(X,\cM)(z):=\sum_{k} {b_{({c})}^{k}}(X,\cM)\cdot (-z)^k \in \bL:=\bZ[z^{\pm 1}],$$
and, respectively, {\it mixed Hodge polynomial} $$h_{({c})}(X,\cM)(y,x,z):=\sum_{p,q,k} {h_{({c})}^{p,q,k}}(X,\cM)\cdot y^px^q(-z)^k  \in \bL:=\bZ[y^{\pm 1},x^{\pm 1},z^{\pm 1}]$$ of $\cM$ in the mixed Hodge module setting. Here, $$b_{({c})}^{k}(X,\cM):=\dim_{\bK} H_{({c})}^k(X,\cM)$$ and $$h_{({c})}^{p,q,k}(X,\cM):=h^{p,q}(H_{({c})}^k(X,\cM)):=\dim_{\bC}Gr^p_F Gr^W_{p+q} H_{({c})}^k(X,\cM)$$ denote the Betti and, respectively, mixed Hodge numbers of the (compactly supported) cohomology $H_{({c})}^*(X,\cM)$ of $\cM$. Similar considerations apply in the context of (\ref{Kut}), where we aim to compute the generating series ${\sum_{n \geq 0} tr_{\Sigma_n}(H^*(X^n)) \cdot t^n}$
in terms of the corresponding Betti numbers and Poincar\'e polynomials of the finite dimensional cohomology $H^*(X)$. 

After composing with the Frobenius character homomorphism \cite{Mc}[Ch.1,Sect.7]:  $$ch_F: C(\Sigma)\otimes_{\bZ}  \bQ:=\bigoplus_n C(\Sigma_n) \otimes_{\bZ}  \bQ \overset{\simeq}{\to} \Lambda \otimes_{\bZ} \bQ = \bQ[p_i, i \geq 1],$$
the generating series 
$$\sum_{n \geq 0} tr_{\Sigma_n}(H^*_{({c})}(X^n,\cM^{\boxtimes n})) \cdot t^n$$
can be regarded as an element in the $\bQ$-algebra $\bL\otimes_{\bZ}  \bQ[p_i, i \geq 1][[t]]$. 
Here, $\Lambda$ is the graded ring of $\bZ$-valued symmetric functions in infinitely many variables $x_m$ ($m \in \bN$), with $p_i=\sum_m x_m^i$ the $i$-th power sum function.

\medskip

The first main result of this note is the following:
\begin{thm}\label{mth} 
For a pair $(X,\cM)$ of a space $X$ with a coefficient $\cM$ as above, the following generating series identity for the Poincar\'e polynomials of characters of external products of $\cM$ holds in the $\bQ$-algebra $\bQ[p_i, i \geq 1,z^{\pm 1}][[t]]$:
\be\label{mf1}
\sum_{n \geq 0} tr_{\Sigma_n}(H^*_{({c})}(X^n,\cM^{\boxtimes n})) \cdot t^n  
 =\exp \left( \sum_{r \geq 1}  p_r \cdot {P_{({c})}}(X,\cM)(z^r) \cdot \frac{t^r}{r}\right).
\ee

Moreover, in the case when $\cM \in D^b\mh(X)$ is a complex of mixed Hodge modules on a complex quasi-projective variety $X$, the following refined 
generating series identity for the mixed Hodge polynomials of characters of external products of $\cM$ holds in the $\bQ$-algebra $\bQ[p_i, i \geq 1, y^{\pm 1},x^{\pm 1},z^{\pm 1}][[t]]$:
\be\label{mf1h}
\sum_{n \geq 0} tr_{\Sigma_n}(H^*_{({c})}(X^n,\cM^{\boxtimes n})) \cdot t^n  
 =\exp \left( \sum_{r \geq 1}  p_r \cdot {h_{({c})}}(X,\cM)(y^r,x^r,z^r) \cdot \frac{t^r}{r}\right).
\ee
\end{thm}

\br A formula similar to (\ref{mf1}) also holds for the characters of the finite dimensional cohomology $H^*(X^n)$ in the context of (\ref{Kut}), i.e., by ``forgetting'' the coefficients in the statement. This also applies to the cases 
$$H^*(X):=H_{({c})}^*(X,\bK) , \ H^*(X,\cO), \ H_{({c})}^*(X,\bQ^H),$$
of cohomology with ``trivial'' coefficients, i.e., the constant sheaf $\bK_X$, structure sheaf $\cO_X$, or the constant mixed Hodge module complex $\bQ_X^H$. (In the last case, one recovers Deligne's mixed Hodge structure on $H_{({c})}^*(X,\bQ)$.)
This fact applies to many results in this paper, and the precise formulation will be left to the interested reader. 
\er

%%%%%%%%%%%%%%%%%%%%%%%%%%%%%%%%%%%

\subsection{Twisting by symmetric group representations}\label{twist}
Additionally, for a fixed $n$, one can consider the coefficient of $t^n$ in the generating series for the characters of cohomology representations $H^*_{({c})}(X^n,\cM^{\boxtimes n})$ of all exterior powers $\cM^{\boxtimes n}$. Moreover, in this case, one can twist the coefficients $\cM^{\boxtimes n}$ by a rational $\Sigma_n$-representation $V$ (see Remark \ref{ptwist}), to get a $\Sigma_n$-equivariant (twisted) coefficient  $V \otimes \cM^{\boxtimes n}$ on $X^n$,  
and compute the corresponding characters for the twisted cohomology $\Sigma_n$-representations $H^*_{({c})}(X^n, V \otimes \cM^{\boxtimes n})$
via the equivariant K\"unneth formula
\be\label{tKue} H^*_{({c})}(X^n, V \otimes \cM^{\boxtimes n})
\simeq V \otimes H^*_{({c})}(X^n,\cM^{\boxtimes n}) \simeq V \otimes H^*_{({c})}(X,\cM)^{\otimes{n}}.\ee
Here, in the mixed Hodge context, we regard $V$ as a pure Hodge structure of type $(0,0)$.
By the multiplicativity of characters, 
we then have:
\be\label{new1} tr_{\Sigma_n}(H^*_{({c})}(X^n, V \otimes \cM^{\boxtimes n}))=tr_{\Sigma_n}(V) \cdot tr_{\Sigma_n}(H^*_{({c})}(X^n, \cM^{\boxtimes n}))\:.\ee
Expanding the exponential series of Theorem \ref{mth}, together with the above multiplicativity, we get our second main result:
\bt\label{tw} In the above notations, the following identity holds in  $\bQ[p_i, i \geq 1,z^{\pm 1}]$:
\be\label{new2a}
tr_{\Sigma_n}(H^*_{({c})}(X^n, V \otimes \cM^{\boxtimes n})) = \sum_{{\lambda=(k_1,k_2, \cdots) \vdash  n} } \frac{p_{\lambda}}{z_{\lambda}} \chi_{\lambda}(V)  \cdot \prod_{r \geq 1} \left( P_{({c})}(H^*(X;\cM)(z^r) \right)^{k_r},
\ee
and similarly for the mixed Hodge context.
Here, for a partition  $\lambda=(k_1,k_2, \cdots)$ of $n$ (i.e., $\sum_{r \geq 1} r \cdot k_r=n$) corresponding to a conjugacy class of an element $\sigma\in \Sigma_n$, we denote by $z_{\lambda}:=\prod_{r \geq 1} r^{k_r} \cdot k_r!$ the order of the stabilizer of $\sigma$, by $\chi_{\lambda}(V)=trace_{\sigma}(V)$ the corresponding trace, and we set $p_{\lambda}:=\prod_{r \geq 1} p_r^{k_r}$. 
\et

\medskip

The formulae of Theorems \ref{mth} and \ref{tw} can be specialized in several different ways, e.g., 
\begin{itemize}
\item[(i)] for specific values of the parameter $z$ (and, resp., $x$, $y$, $z$ in the mixed Hodge context), e.g., the specialization $z=1$ yields Euler-characteristic type formulae; \item[(ii)] for special choices of the coefficient $\cM$ (e.g., intersection cohomology complexes); \item[(iii)] for special values of the Frobenius parameters $p_r$ (e.g., related to symmetric and alternating powers of coefficients);
\item[(iv)] for special choices of the representation $V \in Rep_{\bQ}(\Sigma_n)$, e.g., for $V={\rm Ind}_K^{\Sigma_n}(triv_K)$, the representation induced from the trivial representation of a subgroup $K$ of $\Sigma_n$.
\end{itemize}
These special cases will be discussed in detail in Section \ref{spc}.

%%%%%%%%%%%%%%%%%%%%%%%%%%%%%%%

\subsection{Abstract generating series formulae and Plethysm}

Theorem \ref{mth} is a direct application of a generating series formula for abstract characters $tr_n$ of tensor powers $\cV^{\otimes n}$ of an element $\cV$ in a suitable symmetric monoidal category $(A,\otimes)$, which in our case will be 
$$\cV=H_{({c})}^*(X,\cM), \ \ \ {\rm resp.,} \ \ \ \cV=Gr^*_F Gr^W_*H_{({c})}^*(X,\cM),$$
as an element in the abelian tensor category of finite dimensional (multi-)graded vector spaces. Note that 
 the functor $Gr^*_FGr^W_*$ is an exact tensor functor on the category of mixed Hodge structures, so it is compatible with the K\"unneth isomorphism (\ref{Kue}). 

In more detail, let $A$ be a pseudo-abelian (or Karoubian) $\bQ$-linear additive category which is also symmetric monoidal, with tensor product $\otimes$ $\bQ$-linear additive in both variables. Then the corresponding Grothendieck ring ${\bar K}_0(A)$ of the {\it additive} category $A$ is a pre-lambda ring with a pre-lambda structure defined by (see (\ref{prel2}), and compare also with \cite{Ge, He}):
\begin{equation} \label{prel}
\sigma_t: {\bar K}_0(A)\to {\bar K}_0(A)[[t]]\:,\:\:
[\cV] \mapsto 1+ \sum_{n\geq 1}\; [(\cV^{\otimes n})^{\Sigma_n}] \cdot t^n \:,
\end{equation} 
for $(-)^{\Sigma_n}$ the functor defined by taking the $\Sigma_n$-invariant part (where the pseudo-abelian $\bQ$-linear additive structure of $A$ is used to define the  projector on the invariant part). Recall that a pre-lambda structure on a commutative ring $R$ with unit $1$ is a group homomorphism
$$\sigma_t: (R,+)\to (R[[t]],\cdot)\:;\:\:
r \mapsto 1+ \sum_{n\geq 1}\; \sigma_n(r) \cdot t^n, $$
with $\sigma_1=id_R$, where ``$\cdot$'' on the target side denotes the multiplication of formal power series. 

Let $A_{\Sigma_n}$ be the additive category of the $\Sigma_n$-equivariant objects in $A$, as in \cite{MS}[Sect.4], with corresponding Grothendieck group ${\bar K}_0(A_{\Sigma_n})$. Then one has the following decomposition 
(e.g., see \cite{MS}[Eqn.(45)] and Section \ref{simmon}):
$$
{\bar K}_0(A_{\Sigma_n}) \simeq {\bar K}_0(A)\otimes _{\bZ} Rep_{\bQ}(\Sigma_n), 
$$
with $Rep_{\bQ}(\Sigma_n)$ the ring of rational representations of $\Sigma_n$.
We next denote by $tr_n$ the composition:
$$
\begin{CD}
tr_n:{\bar K}_0(A_{\Sigma_n}) 
\simeq  {\bar K}_0(A)\otimes _{\bZ} Rep_{\bQ}(\Sigma_n) @>id \otimes tr_{\Sigma_n} > >
{\bar K}_0(A)\otimes _{\bZ} C(\Sigma_n).
\end{CD}$$
Fix now an object $\cV \in A$, and consider the generating series:
$$\sum_{n \geq 0} tr_n([\cV^{\otimes n}]) \cdot t^n \in {\bar K}_0(A)  \otimes_{\bZ}  C(\Sigma)[[t]].$$
Let $$cl_n:=(id \otimes ch_F) \circ tr_n,$$
with $ch_F$ denoting the Frobenius character homomorphism as before.

\medskip

In the above notations, the first abstract formula of this paper can now be stated as follows (see Theorem \ref{mta} in Section \ref{simmon}):
\bt\label{mti} For any $\cV \in A$, the following generating series identity holds in the $\bQ$-algebra $\left({\bar K}_0(A) \otimes_{\bZ}  \bQ[p_i, i \geq 1] \right)[[t]] = \left(  \bQ[p_i, i \geq 1]\otimes_{\bZ}  {\bar K}_0(A)\right)[[t]]$:
\be\label{mainabsi} \sum_{n \geq 0} cl_n([\cV^{\otimes n}]) \cdot t^n=\exp \left( \sum_{r \geq 1}  \psi_r([\cV]) \otimes p_r \cdot \frac{t^r}{r}    \right),
\ee
with $\psi_r$ the $r$-th Adams operation of the pre-lambda ring ${\bar K}_0(A)$.  Here we use the Atyiah-type definition of Adams operations given by $ \psi_r([\cV]) :=tr_r([\cV^{\otimes r}])(\sigma_r)$, where $\sigma_r \in \Sigma_r$ is a cycle of length $r$.
\et
Note that by setting $p_r=1$ for all $r$, formula (\ref{mainabsi}) specializes to the well-known pre-lambda ring identity (e.g., see \cite{Kn} or \cite{Mc}[Ch.1,Rem.2.15]):
\be\label{simi} 
\sigma_t\left([\cV]\right)=1+\sum_{n\geq 1}\;  [(\cV^{\otimes n})^{\Sigma_n}] \cdot t^n =\exp\left( \sum_{r\geq 1}\;
\psi_r([\cV]) \cdot \frac{t^r}{r}\right) \in {\bar K}_0(A)\otimes_{\bZ}\bQ[[t]]\:,\ee
relating the pre-lambda structure to the corresponding Adams operations.  Formula (\ref{simi}) was the main tool used for proving our results in \cite{MS} (see also \cite{Ge}). In this paper, we use a more general equivariant approach, which does not rely on the theory of pre-lambda rings.

\medskip

Similarly, Theorem \ref{tw} for twisted coefficients can be derived from the following abstract twisting formula (see Theorem \ref{abstwist} of Sect.\ref{simmon}):
\bt\label{itwist}
For $V$ a rational representation of $\Sigma_n$ and $\cV \in A$, the following identity holds in $\bQ[p_i, i \geq 1] \otimes_{\bZ}  {\bar K}_0(A)$:
\be\label{new3ai}
cl_n(V \otimes \cV^{\otimes n}) = \sum_{{\lambda=(k_1,k_2, \cdots) \vdash n} } \frac{p_{\lambda}}{z_{\lambda}} \chi_{\lambda}(V)  \otimes \prod_{r \geq 1} \left( \psi_r([\cV]) \right)^{k_r},
\ee
where $\chi_{\lambda}(V)=trace_{\sigma}(V)$ for $\sigma \in \Sigma_n$ of cycle-type corresponding to the partition  $\lambda=(k_1,k_2, \cdots)$ of $n$ (i.e., $\sum_{r \geq 1} r \cdot k_r=n$),  $p_{\lambda}:=\prod_{r \geq 1} p_r^{k_r}$ and $z_{\lambda}:=\prod_{r \geq 1} r^{k_r} \cdot k_r!$.
\et

\br In the mixed Hodge context, stronger versions of Theorems \ref{mth} and \ref{tw} are  obtained by regarding $\cV=H^*_{({c})}(X;\cM)$ as an element in the abelian tensor category $A$ of finite dimensional graded rational mixed Hodge structures.
\er

\medskip

In Section \ref{feq}, we indicate further applications of the above abstract setup to suitable equivariant versions of (characters of) Poincar\'e and mixed Hodge polynomials of {\it equivariant} coefficients on a space $X$. This is done by replacing the category $A$ by a suitable category $A_G$ of $G$-equivariant objects in $A$ (as in \cite{MS}), resp., by  the category ${\rm End}(A)$ of endomorphisms of objects in $A$. These categories are still pseudo-abelian $\bQ$-linear tensor categories as before, so that the above abstract formulae apply directly also in this equivariant context.

For instance, we can  consider a compact topological space $X$ with finite dimensional cohomology $H^*(X,\bK)$ (resp., a complex quasi-projective variety $X$), together with the action of:
\begin{enumerate}
\item[(A)] a finite group $G$ acting (algebraically) on $X$, 
\item[(B)] an (algebraic) automorphism $g$  of $X$ of finite order,  
\item[(C)] a (proper algebraic) endomorphism $g:X \to X$. 
\end{enumerate}
Here we only need  the fact that the K\"unneth isomorphism (\ref{Kut}) and (\ref{Kue}) are functorial.

For simplicity, we illustrate here such an equivariant formula 
for the constant coefficient $\bK_X$ on a complex quasi-projective algebraic variety $X$, and for   Macdonald-type generating series for the symmetric products $X^{(n)}$ of $X$ (i.e., with all Frobenius variables $p_r$ set to be equal to $1$), see Theorem \ref{final}: 
 
\bt\label{teq} 
If $g:X \to X$ is a (proper) algebraic endomorphism of a complex quasi-projective variety $X$, then the following equality holds in $\bK[z][[t]]$:
\be\label{eqf}
 \sum_{n \geq 0}  P^g_{({c})}(X^{(n)})(z) \cdot t^n  
=\exp \left( \sum_{r \geq 1}  P^{g^r}_{({c})}(X)(z^r) \cdot \frac{t^r}{r} \right),\ee
with $X^{(n)}:=X^n/{\Sigma_n}$ the $n$-th symmetric product of $X$ and 
$$P^g_{({c})}(X)(z) := \sum_k trace_{g} \left(H^k_{({c})}(X,\bK)  \right) \cdot (-z)^k.$$
A similar formula holds for the equivariant mixed Hodge polynomials.
\et

Note that formula (\ref{eqf}) also holds for an endomorphism $g$ of a compact topological space $X$ with finite dimensional cohomology. 
Formula (\ref{eqf}) specializes for $z=1$ to the usual {\it Lefschetz Zeta function} of the (proper) endomorphism $g: X\to X$.  Moreover, for $g=id_X$ the identity of $X$, formula (\ref{eqf})
reduces to {\it  Macdonald's generating series formula} \cite{Mac} for the Poincar\'e polynomials and Betti numbers of the symmetric products of $X$. For more details, see Section \ref{feq}.

\br 
For the counterpart of the Lefschetz Zeta function
in the context of  constructible sheaves for the Frobenius endomorphism of varieties over finite fields, see also  \cite{SGA5}[Thm. on p.464] and \cite{FK}[Thm.4.4 on p.174]. 
For a similar counterpart of (\ref{eqf}) taking  a weight filtration into account, see \cite{Na}[Prop.8(i)].
 \er
 
\medskip
 
We next apply the above abstract setting in order to develop the plethysm calculus for symmetric sequences with values in a tensor category $(A,\otimes)$ as above. The category of symmetric sequences is defined by $S(A):=\prod_{n \geq 0} A_{\Sigma_n}$, and it is endowed with the (graded) Cauchy product $\odot:={\rm Ind}^{\Sigma_{n+m}}_{\Sigma_n \times \Sigma_m}( \cdot \otimes \cdot)$ given by induction and the monoidal structure.

The category $S(A)$ is pseudo-abelian $\bQ$-linear additive,  with the Cauchy product $\odot$ $\bQ$-linear additive in both variables. 
Moreover, $S(A)$ is also symmetric monoidal with respect to the Cauchy product $\odot$. Then all of our abstract results for $A$ are also applicable in the context of $S(A)$, with 
$${\bar K}_0(S(A)) \simeq \prod_n {\bar K}^{\Sigma_n}_0(A)$$
as ${\bar K}_0(A)$-algebras and commutative graded rings with the induced Cauchy products.  

The graded ring homomorphism 
$$cl:=\sum_n cl_n : \bigoplus_n {\bar K}^{\Sigma_n}_0(A) {\lra} {\bar K}_0(A) \otimes_{\bZ} \bQ[p_i, i \geq 1] $$
(with $p_i$ in degree $i$)  
induces by completion a similar homomorphism
$$cl: {\bar K}_0(S(A)) \lra  {\bar K}_0(A) \otimes_{\bZ} \bQ[[p_i, i \geq 1]].$$

The final abstract formula of this paper is given by the following generalization of a corresponding result of Joyal for the classical case $A={\rm Vect}_{\bQ}$, the category of finite dimensional rational vector spaces (see \cite{J1,J2}).
\bt\label{mt3i} Let $\cV_{\cdot}=(\cV_n)_{n \geq 0}$ , $\cW_{\cdot}=(\cW_n)_{n \geq 0} \in S(A)$ be two symmetric sequences in $A$. Then we have:
\be\label{ma3i}
cl([\cV_{\cdot}] \circ [\cW_{\cdot}])= cl([\cV_{\cdot}] )\circ cl([\cW_{\cdot}]),
\ee
and  
\be\label{ma4i}
cl([\cV_{\cdot}] \ast [\cW_{\cdot}])= cl([\cV_{\cdot}] )\ast cl([\cW_{\cdot}]),
\ee
where in (\ref{ma3i}) one has to assume that $\cW_0=0_A$ is the zero object in $A$.
\et
On the left-hand side of the above formulae (\ref{ma3i}) and (\ref{ma4i}), $\circ$ resp., $\ast$ is the composition pairing (\ref{circus}), resp., the Hadamard product of Remark \ref{Had}(1).
On the right-hand side, $\circ$ resp., $\ast$ denotes the abstract plethysm of Definition \ref{pl}, resp., the internal product of formal power series in the $p_i$'s with coefficients in ${\bar K}_0(A)$, see Remark \ref{Had}(1). 

\medskip

Theorem \ref{mt3i} generalizes similar results of \cite{Ge,GP} for $A$ the category of rational mixed Hodge structures, resp., results of \cite{H}[Sect.4, Sect.5] for $A$ the category of finite-dimensional complex $G$-representations or the category of endomorphisms of finite dimensional ${\overline \bQ}_{\ell}$-vector spaces (applied to the $\ell$-adic cohomology $H^*_c(X,{\overline \bQ}_{\ell})$ of a variety $X$ over a finite field, together with the Frobenius endomorphism). 

\medskip

Let us point out that Theorem \ref{mt3i} is deduced from our two abstract formulae 
(\ref{mainabsi}) and (\ref{new3ai}). On the other hand, as explained in Example \ref{pe} and Remark \ref{Had}(1), these abstract formulae can also be recovered from Theorem \ref{mt3i}.

\medskip

Geometric applications of (\ref{ma3i}) in special cases are given in \cite{Ge,GP} to reduce the calculation of invariants of {\it configuration spaces} $\bF(X,n)$ of $n$ points in $X$ and  their {\it Fulton-MacPherson compactifications} $\bF\bM(X,n)$ (for $X$ smooth) to the calculation of the corresponding invariants for products. This uses the associativity of the composition products 
(see Remark \ref{plinv} for a short discussion, and also  \cite{GP} for some ``motivic'' counterparts). Here we only state the following application to {\it Hilbert schemes of ordered points}, which doesn't seem to be available in the literature (complete details will be explained in a future work by the authors).

Let $X$ be a smooth complex quasi-projective variety of pure dimension $d$. Then the {\it Hilbert scheme $X^{[[n]]}$ of $n$ ordered points in $X$} is defined in \cite{Re}[Def.2.2] as the (reduced) fiber-product via the cartesian diagram
$$\begin{CD} X^{[[n]]} @>>> X^{[n]} \\
@VVV @VV p_n V \\
X^n @>>> X^{(n)}=X^n/\Sigma_n,
\end{CD}$$
with $p_n$ the proper Hilbert-Chow morphism from the usual  Hilbert scheme $X^{[n]}$ of $n$  points in $X$ to the symmetric product $X^{(n)}$. For our application, it is enough to consider $X^{[[n]]}$ (and all other spaces) with its reduced structure so that
$X^{[n]}=X^{[[n]]}/\Sigma_n$ with respect to the induced algebraic $\Sigma_n$-action on $X^{[[n]]}$.
The stratified nature of the Hilbert-Chow morphism discussed in \cite{BBS}[Sect. 2.1] implies the following plethysm formula:

\be\label{hilbord1} \left([H^*_c(X^{[[n]]})]\right)_{n \geq 0}=\left([H^*_c(\bF(X,n))]\right)_{n \geq 0} 
 \circ \Big([H^*_c({\rm Hilb^n_{\mathbb{C}^d,0}})]\Big)_{n \geq 0}.\ee  
Here ${\rm Hilb^n_{\mathbb{C}^d,0}}$ are the corresponding {\it punctual Hilbert schemes} (with 
${\rm Hilb^0_{\mathbb{C}^d,0}}=\emptyset$), viewed with a trivial $\Sigma_n$-action, and 
 $[H^*_c]:=\sum_{i} (-1)^i [H^i_c]\in K_0({\ms})$ in the usual Grothendieck group $K_0({\ms})$ of the Abelian category of mixed Hodge structures.
Moreover, it is important to work with the cohomology with compact support $H^*_c$ and the Euler-characteristic type class 
$[H^*_c]$, since ``additivity'' with respect to decompositions is used in these calculations.
Using the plethistic exponential and  logarithm isomorphism ${\rm Exp}$ and ${\rm Log}$, this can be reformulated 
(by the same argument as in Remark \ref{plinv}) as:

\bex\label{hilbord2} With the above notations, the following equality  holds:
\be\label{hilbord3} 
   \begin{split} \left([H^*_c(X^{[[n]]})]\right)_{n \geq 0} &={\rm Exp} \left( [H^*_c(X)] \otimes {\rm Log}\Big(1+\sum_{i \geq 1} [H^*_c({\rm Hilb^i_{\mathbb{C}^d,0}})]\Big) \right)\\
   &=:\Big(1+\sum_{i \geq 1} [H^*_c({\rm Hilb^i_{\mathbb{C}^d,0}})]\Big)^{[H^*_c(X)]} \in K_0(S({\ms})).
   \end{split}
 \ee \eex

By application of the plethysm formula  (\ref{ma3i}) (or better, its variant for Abelian tensor categories, which is left to the reader),
together with the specialization $p_i\mapsto p_it^i$ for all $i$ (for  specifying the degree), this implies the following equality 
in $\left(K_0({\ms})\otimes_{\bZ} \bQ[p_i, i \geq 1]\right)[[t]]$:
\be \label{hilbord4} 
1+\sum_{n \geq 1} cl_n\left([H^*_c(X^{[[n]]})]\right)t^n = \Big(1+\sum_{i \geq 1} [H^*_c({\rm Hilb^i_{\mathbb{C}^d,0}})]\otimes h_it^i\Big)^{[H^*_c(X)]\otimes 1},
\ee
formulated in terms of the {\it  power structure} in the sense of \cite{GLM2}  associated to the pre-lambda structure on 
$K_0({\ms})\otimes_{\bZ} \bQ[p_i, i \geq 1]$. Here  $h_n \in  \bQ[p_i,  i \geq 1]$ are the
complete symmetric functions ($n\ge 1$).

By application of  the specialization $p_i\mapsto 1$ for all $i$ (corresponding to taking the $\Sigma_n$-invariant part,
with $h_i\mapsto 1$ for all $i$), this implies the following known formula:
\be \label{hilbord5} 
1+\sum_{n \geq 1} [H^*_c(X^{[n]})]t^n = \Big(1+\sum_{i \geq 1} [H^*_c({\rm Hilb^i_{\mathbb{C}^d,0}})]t^i\Big)^{[H^*_c(X)]} \in K_0({\ms})\otimes_{\bZ} \bQ[[t]]
\ee
formulated in terms of the {\it  power structure} associated to the pre-lambda structure on $K_0({\rm mHs})\otimes_{\bZ} \bQ$
(as in \cite{GLM2, CMOSY}, where also the corresponding ``motivic'' formula is discussed). 
Application of the  {\it $E$-polynomial} (as in Section \ref{spc}(i)) recovers the corresponding result of Cheah
\cite{Che}.

%%%%%%%%%%%%%%%%%%%%%%%%%%%%%%%

\subsection{Pseudo-functors}

In the final Section \ref{psf}, we explain the connection of our results  with our previous work from \cite{MS} about generating series of symmetric and alternating powers of suitable coefficients.
In fact, all of this can and will be discussed in the abstract setting of suitable pseudo-functors $A(X)$ with a (derived) pushforward $(-)_*$  as in \cite{MS}.  For simplicity, we focus here only  on the complex quasi-projective context, and the following examples:

\begin{enumerate}
\item[(a)] $A(X)=D^b_{f}(X;\bK)$ is the triangulated category consisting of bounded complexes of sheaves of $\bK$-vector spaces with finite dimensional compactly supported cohomology, viewed as a pseudo-functor with respect to $(-)_!$.
\item[(b)] $A(X)=D^b_c(X;\bK)$ is the derived category of bounded complexes of sheaves of $\bK$-vector spaces with constructible cohomology,
viewed as a pseudo-functor with respect to either of the pushforwards $(-)_*$ or $(-)_!$.
\item[(c)] $A(X)=D^b_{coh}(X)$ is the derived category of bounded complexes of sheaves of $\cO_X$-modules with coherent cohomology, with $X$ projective, and
viewed as a pseudo-functor with respect to   $(-)_*=(-)_!$.
\item[(d)]$A(X)=D^b\mh(X)$ is the bounded derived category of algebraic mixed Hodge modules, viewed as a pseudo-functor with respect to either of the pushforwards $(-)_*$ or $(-)_!$.
\end{enumerate}

Note that in all these cases $A(X)$ is a pseudo-abelian (or Karoubian) $\bQ$-linear additive category, such that $A(pt)$ 
(for $X=pt$ a point space) is endowed with a $\bQ$-linear tensor structure $\otimes$, which makes it into a symmetric monoidal category as before.  

\medskip

So we can apply our results also to the (derived) pushforwards $\cV:=k_*\cM \in A(pt)=:A$ for $\cM\in A(X)$,
with $k$ a constant map to the point $pt$. By the corresponding equivariant derived  K\"{u}nneth formula, these results can also be formulated in terms of suitable  $\Sigma_n$-equivariant external products $\cM^{\boxtimes n}\in A(X^n)$ of $\cM\in A(X)$.
Finally, one gets  symmetric and alternating powers $\cM^{(n)}, \cM^{\{n\}}\in A(X^{(n)})$, as well as Schur objects
$S_V(\cM)\in A(X^{(n)})$, as coefficients on the symmetric product $X^{(n)}=X^n/\Sigma_n$ in such a way that
$$(k_*\cM)^{(n)}=k_*(\cM^{(n)}), \ (k_*\cM)^{\{n\}}=k_*(\cM^{\{n\}}) \quad \text{and} \quad S_V(k_*\cM)=k_*(S_V\cM)\:.$$

These induced ($\Sigma_n$-equivariant) coefficients on $X^{(n)}$ (resp., $X^n$) become important in our related work on
(equivariant) characteristic classes as in \cite{CMOSY, CMSSY, MS15}. For $X$ projective, one can recover  some of the
Euler-characteristic type  results of this paper by taking the degree of similar formulae for these  characteristic classes.  Note that in the mixed Hodge context, these characteristic class formulae only take into account the Hodge filtration.
So instead of the mixed Hodge polynomial $h(X,\cM)(y,x,z)$, one can only recover at the degree level the corresponding $\chi_{y}$-genus
$$\chi_{y}(X,\cM):=h(X,\cM)(y,x=1,z=1)\:.$$
For $X$  a complex projective manifold and $\cM=\bQ_X^H$ the constant (shifted)  Hodge module, this
becomes the classical {\it Hirzebruch $\chi_y$-genus}. This is also the reason why we label by $y$ the parameter corresponding to the Hodge filtration (hence the unusual ordering $y,x,z$ of parameters in the definition of the mixed Hodge polynomial).

\medskip

Finally, note that there is of course some overlap of our paper with the unpublished preprint \cite{Ge}, which uses a (pre-) lambda ring approach and focuses in the applications mainly on the mixed Hodge context and Euler-characteristic type invariants.
Instead, in this paper we use an equivariant appoach, relying only on the Atiyah-type description of the Adams operations 
(which fits nicely with the classical Frobenius character homomorphism), together with the use of Schur functors.
This equivariant point of view is motivated by, and closely related to, the corresponding deeper results on (equivariant) characteristic classes, 
as in our papers \cite{CMOSY, CMSSY, MS15} (but which are not needed in the work presented here). Moreover, we take as a starting point our two abstract Theorems \ref{mti} and \ref{itwist}, and deduce everything else form these, including the abstract  plethysm results  in Theorem \ref{mt3i}
(which provide a generalization of the corresponding results of Joyal  \cite{J1,J2}). We also apply our basic abstract results to many different concrete situations, by working with different coefficients and in various situations, also in a derived or equivariant context. All of these applications, including the
abstract plethysm results, are not contained in \cite{Ge}.

\begin{ack} 
L. Maxim was partially supported by grants from NSF, NSA, by a fellowship from the Max-Planck-Institut f\"ur Mathematik,  Bonn, and by the Romanian Ministry of National Education, CNCS-UEFISCDI, grant PN-II-ID-PCE-2012-4-0156.
J. Sch\"urmann was supported by the SFB 878 ``groups, geometry and actions". 
\end{ack}

%%%%%%%%%%%%%%%%%%%%%%%%%%%%%%%%% 
%%%%%%%%%%%%%%%%%%%%%%%%%%%%%%%%% 

\section{Abstract generating series identities and Plethysm}\label{abssec}
In this section, we prove the abstract Theorems \ref{mti} and \ref{itwist}
 from the Introduction. Moreover, we develop an abstract plethysm calculus for symmetric sequences in a suitable pseudo-abelian $\bQ$-linear additive tensor category $A$. 
 
\subsection{Symmetric monoidal categories and abstract generating series}\label{simmon}
Let $A$ be a pseudo-abelian (or Karoubian) $\bQ$-linear additive category which is also symmetric monoidal, with the tensor product $\otimes$ $\bQ$-linear additive in both variables. Let $${\bar K}_0(A):=\big(K_0(A),\oplus \big)$$ denote the corresponding Grothendieck ring of the {\it additive} category $A$. 
Similarly, let $A_{\Sigma_n}$ be the additive category of the $\Sigma_n$-equivariant objects in $A$, as in \cite{MS}[Sect.4], with corresponding Grothendieck group $${\bar K}^{\Sigma_n}_0(A):={\bar K}_0(A_{\Sigma_n}).$$ Then one has the following decomposition (e.g., see \cite{MS}[Eqn.(45)]):
\begin{equation}\label{dec}
{\bar K}^{\Sigma_n}_0(A) \simeq {\bar K}_0(A)\otimes _{\bZ} Rep_{\bQ}(\Sigma_n) \simeq Rep_{\bQ}(\Sigma_n) \otimes _{\bZ}  {\bar K}_0(A), 
\end{equation}
with $Rep_{\bQ}(\Sigma_n)$ the ring of rational representations of $\Sigma_n$.
In fact, this follows directly from the corresponding decomposition of
$\cY \in A_{\Sigma_n}$ by  {\it Schur functors} 
$S_{\mu}:A_{\Sigma_n}\to A , \  \cY \mapsto (V_{\mu} \otimes \cY)^{\Sigma_n}$
(e.g., see \cite{De,He}):
\be\label{Sch} \cY \simeq \sum_{\mu \vdash n} \; V_{\mu}\otimes S_{\mu}(\cY)\:,\ee
with $V_{\mu}\simeq V^*_{\mu} $ the (self-dual) irreducible $\bQ$-representation of $\Sigma_n$ corresponding to the partition $\mu$ of $n$. Here, the Karoubian $\bQ$-linear additive structure of $A$ is used to defined the $\Sigma_n$-invariant part functor by the projector $$(-)^{\Sigma_n}:=\frac{1}{n!} \sum_{\sigma \in \Sigma_n} \sigma_*,$$ with $\sigma_*$ denoting the action of $\sigma \in \Sigma_n$.

As in the classical representation theory, the rings ${\bar K}^{\Sigma_n}_0(A)$ have product, induction and restriction functors compatible with (\ref{dec}), which are induced from the corresponding functors on $A_{\Sigma_n}$, see \cite{De}[Sect.1], \cite{He}[Sect.4.1]:
\begin{enumerate}
\item[(a)] the product: $$ \otimes: {\bar K}^{\Sigma_n}_0(A) \otimes_{\bZ} {\bar K}^{\Sigma_m}_0(A) \to {\bar K}^{\Sigma_{n} \times \Sigma_{m}}_0(A)$$
induced from
$$\otimes: A_{\Sigma_n}\otimes A_{\Sigma_m} \to A_{\Sigma_{n} \times \Sigma_{m}}.$$
\item[(b)] induction functor: $${\rm Ind}_{\Sigma_n \times \Sigma_m}^{\Sigma_{n+m}}: {\bar K}^{\Sigma_{n} \times \Sigma_{m}}_0(A) \to {\bar K}^{\Sigma_{n+m}}_0(A)$$
induced from the additive functor 
$${\rm Ind}_{\Sigma_n \times \Sigma_m}^{\Sigma_{n+m}}: A_{\Sigma_n} \times A_{\Sigma_m} \to A_{\Sigma_{n+m}}, 
\ \ \cY \mapsto (\bQ[\Sigma_{n+m}] \otimes \cY)^{\Sigma_{n} \times \Sigma_m}.$$
\item[(c)]  the restriction functor $${\rm Res} _{\Sigma_n \times \Sigma_m}^{\Sigma_{n+m}}: {\bar K}_0^{\Sigma_{n+m}}(A) \to {\bar K}_0^{\Sigma_{n} \times \Sigma_{m}}(A)$$
induced from the obvious restriction functor: ${\rm Res} _{\Sigma_n \times \Sigma_m}^{\Sigma_{n+m}}: A_{\Sigma_{n+m}} \to A_{\Sigma_{n} \times \Sigma_{m}}$.
Similarly, for the restriction functor
$${\rm Res} ^{\Sigma_n \times \Sigma_n}_{\Sigma_{n}}: {\bar K}_0^{\Sigma_{n} \times \Sigma_{n}}(A) \to {\bar K}_0^{\Sigma_{n}}(A)$$ induced from the diagonal inclusion $\Sigma_n \hookrightarrow \Sigma_{n} \times \Sigma_{n}$. 

\end{enumerate} 

\medskip

We denote by $tr_n$ the composition:
$$
\begin{CD}
tr_n:{\bar K}^{\Sigma_n}_0(A) 
\simeq  {\bar K}_0(A)\otimes _{\bZ} Rep_{\bQ}(\Sigma_n) @>id \otimes tr_{\Sigma_n} > >
{\bar K}_0(A)\otimes _{\bZ} C(\Sigma_n),
\end{CD}$$
where $tr_{\Sigma_n}$ associates to a representation its character.
By the above considerations, $tr_n$ is compatible with the product, induction and restriction functors, with the corresponding classical notions for the character group. Therefore, we get an induced graded ring homomorphism (which becomes an isomorphism after tensoring with $\bQ$)   $$tr:=\sum_n tr_n : \bigoplus_n {\bar K}^{\Sigma_n}_0(A) {\lra}  {\bar K}_0(A)\otimes _{\bZ}  \left( \bigoplus_n C(\Sigma_n) \right)= {\bar K}_0(A)\otimes _{\bZ} C(\Sigma).  $$ Here the  commutative induction product on both sides is given by: $$\odot:={\rm Ind}^{\Sigma_{n+m}}_{\Sigma_n \times \Sigma_m}(\cdot \otimes \cdot).$$

\medskip

Fix now an object $\cV \in A$, and consider the generating series:
$$\sum_{n \geq 0} tr_n([\cV^{\otimes n}]) \cdot t^n \in {\bar K}_0(A)  \otimes_{\bZ} C(\Sigma)[[t]].$$ 
\br\label{wd} Note that the total power maps $$\cV \mapsto \sum_{n \geq 0} [\cV^{\otimes n}] \cdot t^n \mapsto \sum_{n \geq 0}tr_n( [\cV^{\otimes n}]) \cdot t^n$$ only depend on the Grothendieck class $[\cV] \in {\bar K}_0(A)$. This follows as in \cite{MS}[Prop.3.2] from the identity (see \cite{De})
\be\label{D} 
(\cV \oplus \cV')^{\otimes n} \simeq \bigoplus_{i+j=n} {\rm Ind}_{\Sigma_i \times \Sigma_j}^{\Sigma_{n}} \left( \cV^{\otimes i} \otimes {\cV'}^{\otimes j} \right),
\ee for $\cV, \cV' \in A$.
\er

Let $$cl_n:=(id \otimes ch_F) \circ tr_n$$ and \be\label{cln} cl:=\sum_n cl_n : \bigoplus_n {\bar K}^{\Sigma_n}_0(A) {\lra} {\bar K}_0(A) \otimes_{\bZ} \bQ[p_i, i \geq 1] \ee be the transformations obtained 
by composing $tr_n$ and, resp,. $tr$ with the Frobenius character homomorphism  $$ch_F: C(\Sigma)\otimes_{\bZ} \bQ=\bigoplus_n C(\Sigma_n) \otimes_{\bZ} \bQ \overset{\simeq}{\lra} \bQ[p_i, i \geq 1].$$
Note that $cl$ is a graded ring homomorphism.
Then 
 the generating series $\sum_{n \geq 0} cl_n([\cV^{\otimes n}])  \cdot t^n$ is an element in the formal power series ring of the $\bQ$-algebra ${\bar K}_0(A) \otimes_{\bZ} \bQ[p_i, i \geq 1].$ Note that the homomorphisms $$
  {\bar K}_0(A)  \otimes_{\bZ} \left(\bigoplus_n Rep_{\bQ}(\Sigma_n)\right)[[t]]  \to
 {\bar K}_0(A)  \otimes_{\bZ} C(\Sigma)[[t]]  \to {\bar K}_0(A)  \otimes_{\bZ}   C(\Sigma)  \otimes_{\bZ} \bQ[[t]] $$ are injective if ${\bar K}_0(A)$ is $\bZ$-torsion-free, so no information is lost in this case after tensoring with $\bQ$ or after applying the Frobenius character homomorphism, hence $cl$ is injective in this situation. For example, this is the case if $A$ is the tensor category of finite dimensional multi-graded vector spaces, or the  category of (polarizable) mixed Hodge structures.
 
\medskip

We can now state our first abstract generating series formula:
\bt\label{mta} For any $\cV \in A$, the following generating series identity holds in the $\bQ$-algebra 
$\left({\bar K}_0(A) \otimes_{\bZ} \bQ[p_i, i \geq 1] \right)[[t]] = \left(\bQ[p_i, i \geq 1]\otimes_{\bZ} {\bar K}_0(A)    \right)[[t]] $:
\be\label{mainabsii} \sum_{n \geq 0}  cl_n([\cV^{\otimes n}])  \cdot t^n=\exp \left( \sum_{r \geq 1}  \psi_r([\cV]) \otimes p_r \cdot \frac{t^r}{r}    \right),
\ee
with $\psi_r$ the $r$-th Adams operation of the pre-lamda ring ${\bar K}_0(A)$. Here we use the Atyiah-type definition of Adams operations given by $ \psi_r([\cV]) :=tr_r([\cV^{\otimes r}])(\sigma_r)$, where $\sigma_r \in \Sigma_r$ is a cycle of length $r$. 
\et

\begin{proof}
This formula can be seen as a special case of Theorem 3.1 from our previous work \cite{MS15}. However,  here we give a direct proof based on the calculus of symmetric functions, adapted to the context of this section.
 
For $\sigma \in \Sigma_n$, we denote by $$tr_n([\cV^{\otimes n}])(\sigma) \in {\bar K}_0(A) $$ the element obtained from $tr_n([\cV^{\otimes n}])$ by evaluating the character at $\sigma$. Then, if $\sigma \in \Sigma_n$ has cycle-type $(k_1,k_2, \cdots )$, by using the fact that $tr_n$ commutes with the restriction and product functors 
it follows that the following multiplicativity property holds:
\be\label{multp}tr_n([\cV^{\otimes n}])(\sigma)=\bigotimes_{r}\left( tr_r([\cV^{\otimes r}])(\sigma_r) \right)^{k_r},\ee
where $\sigma_r \in \Sigma_r$ is a cycle of length $r$.
For any $r \geq 1$, let us now set $$b_r:=tr_r([\cV^{\otimes r}])(\sigma_r) \in {\bar K}_0(A).$$ 
By the definition of the Frobenius character \cite{Mc}[Ch.1, Sect.7], we have:
\be\label{defF}
 cl_n([\cV^{\otimes n}]) =\frac{1}{n!} \sum_{\sigma \in \Sigma_n} tr_n([\cV^{\otimes n}])(\sigma) \otimes \psi(\sigma) \in {\bar K}_0(A) \otimes_{\bZ} \bQ[p_i, i \geq 1],
\ee
where $$\psi(\sigma)=\prod_{r \geq 1} p_r^{k_r}=p_{\lambda}$$ for $\sigma \in \Sigma_n$ in the conjugacy class corresponding to the partition $\lambda:=(k_1,k_2, \cdots )$ of $n$ (i.e., $\sum_{r \geq 1} rk_r=n$). 
Then by (\ref{multp}), formula  (\ref{defF}) can be re-written as:
\be\label{defF1}
 cl_n([\cV^{\otimes n}]) = \sum_{\lambda \vdash  n} \frac{p_\lambda}{z_{\lambda} }\otimes \prod_{r \geq 1}  b_r ^{k_r} \in {\bar K}_0(A) \otimes_{\bZ} \bQ[p_i, i \geq 1],
\ee
with $z_{\lambda}:=\prod_{r \geq 1} r^{k_r} \cdot k_r !$ the order of the stabilizer in $\Sigma_n$ of an element of cycle-type $\lambda$. 
So, we have as in \cite{Mc}[p.25] (see also \cite{Ka}[p.554]):
\be
\begin{split}
\exp \left( \sum_{r \geq 1}  b_r \otimes p_r \cdot \frac{t^r}{r}    \right) 
&= \prod_{r \geq 1} \exp \left( b_r \otimes p_r \cdot \frac{t^r}{r}   \right) \\
&= \prod_{r \geq 1} \sum_{k_r=0}^{\infty} \frac{(b_r \otimes p_r)^{k_r}}{ r^{k_r} \cdot k_r!} \cdot t^{rk_r}\\
&=\sum_{n \geq 0} \left( \sum_{\lambda \vdash n}  \frac{p_{\lambda}}{z_\lambda} \otimes \prod_{r \geq 1}  b_r ^{k_r} \right) \cdot t^n \\
&\overset{(\ref{defF1})}{=} \sum_{n \geq 0} cl_n([\cV^{\otimes n}])  \cdot t^n.
\end{split}
\ee
To conclude the proof of the theorem, recall from \cite{MS}[Sect.3] that the $r$-th Adams operation on ${\bar K}_0(A)$ can be given as $$\psi_r([\cV])=tr_r([\cV^{\otimes r}])(\sigma_r)=:b_r,$$ for $\sigma_r$ a cycle of length $r$ in $\Sigma_r$ (as originally introduced by Atiyah in the context of topological $K$-theory \cite{At}).
\end{proof}

\medskip

We next explain in this abstract setting the twisting construction used in the Introduction (see Section \ref{twist}). Let ${\rm Vect}_{\bQ}(\Sigma_n)$ be the category of  finite dimensional rational $\Sigma_n$-representations. We define a pairing 
\be\label{apair} {\rm Vect}_{\bQ}(\Sigma_n) \times A_{\Sigma_n} \overset{\otimes}{\lra} A_{\Sigma_n}; \ (V, \cY) \mapsto V \otimes \cY\ee
by the composition 
$${\rm Vect}_{\bQ}(\Sigma_n) \times A_{\Sigma_n} \overset{\otimes}{\lra} A_{\Sigma_n \times \Sigma_n}\overset{{\rm Res}}{\lra} A_{\Sigma_n} ,$$
with the underlying tensor product $\otimes$ defined via the $\bQ$-linear additive structure of $A$ (as in \cite{De}) together with its induced $\Sigma_n$-action on each factor, and ${\rm Res}:={\rm Res}^{\Sigma_n \times \Sigma_n}_{\Sigma_n}$ denoting the restriction functor for the diagonal subgroup $\Sigma_n \hookrightarrow \Sigma_n \times \Sigma_n$. This induces a pairing 
\be\label{pair} Rep_{\bQ}(\Sigma_n) \times {\bar K}^{\Sigma_n}_0(A) \overset{\otimes}{\lra} {\bar K}^{\Sigma_n}_0(A)\ee
on the corresponding Grothendieck groups such that 
\be\label{npairing} tr_n([V \otimes \cY])=tr_{\Sigma_n}(V) \cdot tr_n([\cY]) \in {\bar K}_0(A) \otimes_{\bZ} C(\Sigma_n) \simeq C(\Sigma_n)  \otimes_{\bZ}  {\bar K}_0(A),\ee
for $V$ a rational $\Sigma_n$-representation and $\cY \in A_{\Sigma_n}$, with multiplication $\cdot$ induced by the usual multiplication of class functions.

By using formula (\ref{defF1}), together with the above multiplicativity (\ref{npairing}), we obtain (after composing with the Frobenius character homomorphism $ch_F$) our second abstract result:
\bt\label{abstwist}
In the above notations, the following identity holds in $\bQ[p_i, i \geq 1] \otimes_{\bZ }{\bar K}_0(A)$:
\be\label{new3a}
cl_n(V \otimes \cV^{\otimes n}) = \sum_{{\lambda=(k_1,k_2, \cdots) \vdash  n} } \frac{p_{\lambda}}{z_{\lambda}} \chi_{\lambda}(V)  \otimes \prod_{r \geq 1} \left( \psi_r([\cV]) \right)^{k_r},
\ee
where $\chi_{\lambda}(V)=trace_{\sigma}(V)$ for $\sigma \in \Sigma_n$ of cycle-type corresponding to the partition  $\lambda=(k_1,k_2, \cdots)$ of $n$ (i.e., $\sum_{r \geq 1} r \cdot k_r=n$),  $p_{\lambda}:=\prod_{r \geq 1} p_r^{k_r}$,  $z_{\lambda}:=\prod_{r \geq 1} r^{k_r} \cdot k_r!$, and $\psi_r$ the $r$-th Adams operation on ${\bar K}_0(A)$ as before.
\et

\bex\label{e1n} If $A={\rm Vect}_{\bQ}$ is the category of finite dimensional rational vector spaces, then for $\cV=\bQ$ the unit of $A$ with respect to the tensor product the above formula (\ref{new3a}) specializes to the following well-known description of the homomorphism
$$cl_n:Rep_{\bQ}(\Sigma_n) \lra \bQ[p_i, i \geq 1]$$
given by (see also \cite{Mc}[Sect.7, (7.2)]):
 \be\label{e2n}
cl_n(V) = \sum_{{\lambda=(k_1,k_2, \cdots) \vdash  n} } \frac{\chi_{\lambda}(V)}{z_{\lambda}} \cdot \prod_{r \geq 1} p_r^{k_r}.
\ee
\eex

We next make the following definition (compare \cite{Mac1}).
\bd\label{defpow} Let $V$ be a finite dimensional rational $\Sigma_n$-representation. The associated {\it Schur} (or {\it homogeneous polynomial}) {\it functor} $S_V:A \to A$ is defined by 
$$S_V(\cV):=(V \otimes \cV^{\otimes n})^{\Sigma_n}.$$
If $V=V_{\mu} \simeq V^*_{\mu}$ is the (self-dual)  irreducible representation of $\Sigma_n$ corresponding to a partition $\mu$ of $n$, we denote by $S_{\mu}:=S_{V_{\mu}}$ the corresponding Schur functor. If $V$ is the trivial, resp., the sign $\Sigma_n$-representation, the corresponding Schur functor is the {\it $n$-th symmetric power}, resp., the {\it $n$-th alternating power}  of $\cV$, denoted by $\cV^{(n)}$ and, resp.,  $\cV^{\{n\}}$. 
\ed

\bex\label{indfr} Let $K$ be a subgroup of $\Sigma_n$, and $W$ a $K$-representation. 
The Schur functor associated to the induced $\Sigma_n$-representation $V:={\rm Ind}^{\Sigma_n}_K(W)$ is given (by using Frobenius reciprocity) by:
$$S_V(\cV):=(W \otimes \cV^{\otimes n})^{K},$$
where $\cV^{\otimes n}$ is considered with the restricted $K$-action.
\eex

\br 
The Schur functor $S_V$ associated to $V$ induces a corresponding pairing (additive  in the first factor)
$$\circ: Rep_{\bQ}(\Sigma_n) \times {\bar K}_0(A) \lra {\bar K}_0(A)$$
on Grothendieck groups, defined via the composition:
$$Rep_{\bQ}(\Sigma_n) \times {\bar K}_0(A) \overset{id \times (-)^{\otimes n}}{\lra} Rep_{\bQ}(\Sigma_n) \times {\bar K}^{\Sigma_n}_0(A) \overset{\otimes}{\lra} {\bar K}^{\Sigma_n}_0(A) \overset{(-)^{\Sigma_n}}{\lra} {\bar K}_0(A),$$
where the $n$-th power map
$$(-)^{\otimes n} :{\bar K}_0(A) \to {\bar K}^{\Sigma_n}_0(A); \ [\cV] \mapsto [\cV^{\otimes n}]$$
is well-defined by Remark \ref{wd}, $\otimes$ is the pairing defined above, and ${\bar K}^{\Sigma_n}_0(A) \overset{(-)^{\Sigma_n}}{\lra} {\bar K}_0(A)$ is induced from the corresponding additive projection functor.
\er

By specializing (\ref{new3a}) to $p_r=1$ for all $r$ (which, by the Schur functor decomposition (\ref{Sch}), corresponds to taking the $\Sigma_n$-invariant part), we obtain a computation of the Grothendieck class $[S_V(\cV)]=S_V([\cV])$ of the Schur  functor associated to $V$ in terms of Adams operations. More precisely, 
\bc In the above notations, we have:
\be\label{new3b}
S_V([\cV])= \sum_{{\lambda=(k_1,k_2, \cdots) \vdash  n} } \frac{1}{z_{\lambda}} \chi_{\lambda}(V)  \otimes \prod_{r \geq 1} \left( \psi_r([\cV]) \right)^{k_r} \in {\bar K}_0(A) \otimes_{\bZ} \bQ.
\ee 
\ec

\br\label{r3n} In view of Example \ref{e1n}, formula (\ref{new3b}) can be re-written by using the substitution homomorphism $$\bQ[p_i, i \geq 1] \lra  {\bar K}_0(A) \otimes_{\bZ} \bQ \  , \ \ \ p_i \mapsto \psi_i([\cV])$$ as follows:
\be\label{sub1}
\begin{split}
S_V([\cV]) &=cl_n(V)(\psi_1([\cV]), \cdots, \psi_r([\cV]), \cdots) \\
&= \left( \sum_{{\lambda=(k_1,k_2, \cdots) \vdash n} } \frac{1}{z_{\lambda}} \chi_{\lambda}(V)  \cdot \prod_{r \geq 1} p_r^{k_r} \right) \Bigg\rvert_{p_r \mapsto \psi_r([\cV])}.
\end{split}
\ee 
\er

Finally, the Schur functor decomposition (\ref{Sch}) yields (after composing with the Frobenius character homomorphism $ch_F$), the following identity for any $\cV \in A$:
\be\label{new3c}
cl_n([\cV^{\otimes n}])=\sum_{{\mu \vdash n} } s_{\mu} \otimes S_{\mu}([\cV]) \in \Lambda \otimes_{\bZ} {\bar K}_0(A),
\ee
with $s_{\mu}:=ch_F(V_{\mu}) \in \Lambda \subset \bQ[p_i, i \geq 1]$ the corresponding {\it Schur functions}, see \cite{Mc}[Ch.1, Sect.3 and Sect.7]. Note that the Frobenius character $ch_F$ induces an isomorphism of graded rings 
$$ch_F : Rep_{\bQ}(\Sigma):= \bigoplus_n Rep_{\bQ}(\Sigma_n) \overset{\simeq}{\lra} \Lambda \subset \bQ[p_i, i \geq 1].$$

\br\label{corr} The non-degenerate pairing $Rep_{\bQ}(\Sigma_n) \times Rep_{\bQ}(\Sigma_n) \lra \bZ$, given by $(V,W)\mapsto \dim_{\bQ}(V \otimes W^*)^{\Sigma_n}$ induces a duality isomorphism 
$$D:Rep_{\bQ}(\Sigma_n) \simeq {\rm Hom}_{\bZ}(Rep_{\bQ}(\Sigma_n),\bZ)=:Rep_{\bQ}(\Sigma_n)_*$$ identifying the Schur functor $S_V:{\bar K}_0(A) \to {\bar K}_0(A)$ with the corresponding {\it operation} 
on ${\bar K}_0(A)$ defined by $D(V)$, as in \cite{MS}[Sect.3] (where we followed Atiyah's approach to K-theory operations). Here, $W^*$ denotes the dual representation of $W$. 
Summing over all $n$, we get isomorphisms of commutative graded rings 
$$\begin{CD}\Lambda @< ch_F< \sim < Rep_{\bQ}(\Sigma) @> D > \sim > Rep_{\bQ}(\Sigma)_* \end{CD}$$
identifying their respective operations on ${\bar K}_0(A)$ (see also \cite{Bi2}[Lem.2.6] and \cite{UDV}[Cor.5.2]). Here, 
\begin{enumerate}
\item $\Lambda$ acts as a universal lambda ring on ${\bar K}_0(A)$, as in \cite{Ge}.
\item $Rep_{\bQ}(\Sigma)$ acts via direct sums of Schur functors (also called polynomial functors in \cite{Mac1}), as considered in the present paper, i.e., with pairing
\be\label{circ} \circ: Rep_{\bQ}(\Sigma) \times {\bar K}_0(A) \lra {\bar K}_0(A).\ee
Moreover, for $[\cV] \in {\bar K}_0(A)$, the induced evaluation map $-\circ [\cV]:Rep_{\bQ}(\Sigma)  \to  {\bar K}_0(A)$ is (by Frobenius reciprocity) a {\it ring homomorphism}.
\item $Rep_{\bQ}(\Sigma)_*$ acts via operations as in \cite{MS}[Sect.3].
\end{enumerate}
\er

\bex Let us illustrate the above correspondence from Remark \ref{corr} in some important situations. 
\begin{itemize}
\item[(a)] The trivial $\Sigma_n$-representation $triv_n \in Rep_{\bQ}(\Sigma_n)$ corresponds by Example \ref{e1n} to the complete symmetric function $h_n:=cl_n(triv_n) \in \Lambda \subset \bQ[p_i,  i \geq 1]$, respectively, to the homomorphism $\sigma_n:=\frac{1}{n!}\sum_{\sigma \in \Sigma_n} trace_{\sigma} \in Rep_{\bQ}(\Sigma_n)_*$ (as in \cite{MS}[Sect.3]).
\item[(b)] The sign $\Sigma_n$-representation $sign_n \in Rep_{\bQ}(\Sigma_n)$ corresponds by Example \ref{e1n} to the elementary symmetric function $e_n:=cl_n(sign_n) \in \Lambda \subset \bQ[p_i,  i \geq 1]$, respectively, to the homomorphism $\lambda_n:=\frac{1}{n!}\sum_{\sigma \in \Sigma_n} (-1)^{sign(\sigma)} \cdot trace_{\sigma} \in Rep_{\bQ}(\Sigma_n)_*$ (as in \cite{MS}[Sect.3]).
\item[({c})] Let $\lambda=(k_1,k_2,\cdots) \vdash n$ be a partition of $n$. The representation ${\rm Ind}^{\Sigma_n}_{\prod_r (\Sigma_r)^{k_r}} (triv)$ induced from the tensor product of trivial representations corresponds to the symmetric function $h_{\lambda}:=\prod_r h_r^{k_r}$. Similarly, the representation ${\rm Ind}^{\Sigma_n}_{\prod_r (\Sigma_r)^{k_r}} (sign)$ induced from the tensor product of the sign representations corresponds to the symmetric function $e_{\lambda}:=\prod_r e_r^{k_r}$.
\end{itemize}
\eex

In this paper, we do not make use of (1) from the list in Remark \ref{corr}, as it was done in \cite{Ge}. Instead, we take an equivariant approach (as in \cite{MS15}), relying only on the description of Adams operations fitting with (3). The use of Schur functors will be essential in the next section, for the study of the abstract plethysm calculus in a pseudo-abelian $\bQ$-linear tensor category $A$ as before.  Moreover, Frobenius reciprocity and (\ref{D}) directly imply that 
 the Grothendieck ring ${\bar K}_0(A)$ is a pre-lambda ring,  
 with pre-lambda structures defined by 
\begin{equation} \label{prel2}
\sigma_t: {\bar K}_0(A)\to {\bar K}_0(A)[[t]]\:,\:\:
[\cV] \mapsto 1+ \sum_{n\geq 1}\; [\cV^{(n)}] \cdot t^n \:,
\end{equation} 
respectively,
\begin{equation} \label{prel3}
\lambda_t: {\bar K}_0(A)\to {\bar K}_0(A)[[t]]\:,\:\:
[\cV] \mapsto 1+ \sum_{n\geq 1}\; [\cV^{\{n\}}] \cdot t^n \:.
\end{equation} 

\br Our notion of Adams operations agrees by \cite{At}[Corollary 1.8] with the Adams operations of the above pre-lambda structures. These pre-lambda structures on  ${\bar K}_0(A)$ are in fact opposite lambda ring structures  (see \cite{Ge,He}). However, the latter property will not be used in this paper. 
\er

As an illustration of our equivariant approach, we give here a direct proof of the following result (which can also be obtained from the above mentioned lambda ring structure of ${\bar K}_0 (A)$):
\bl\label{l1n} The $r$-th Adams operation $\psi_r:{\bar K}_0 (A) \to {\bar K}_0 (A)$ is a ring homomorphism, for any positive integer $r$, with $\psi_1$ the identity.
\el
\begin{proof}
Let $1_A$ be the unit of $A$ with respect to $\otimes$. Then, $\psi_r([1_A])=[1_A]$ since 
$$1_A^{\otimes n} \simeq 1_A \otimes \bQ^{\otimes n}.$$
Similarly, $\psi_r([\cV] \otimes [\cV'])=\psi_r([\cV]) \otimes \psi_r([\cV'])$ follows by the multiplicativity of $tr_n$ together with the $\Sigma_n$-equivariant isomorphism $$(\cV \otimes \cV')^{\otimes n} \simeq \cV^{\otimes n}  \otimes \cV'^{\otimes n}.$$
Finally, the identity $\psi_r([\cV] \oplus [\cV'])=\psi_r([\cV]) \oplus \psi_r([\cV'])$ follows from (\ref{D}) together with  $$tr_r\left({\rm Ind}_{\Sigma_i \times \Sigma_j}^{\Sigma_{r}} \left( \cV^{\otimes i} \otimes {\cV'}^{\otimes j} \right)\right)(\sigma_r)=0,$$ whenever $i$ or $j$ are non-zero. 
\end{proof}

\br The proof of Lemma \ref{l1n} is an abstract adaptation of Atiyah's argument given in the special case when $A$ is the category of vector bundles on a topological space, see \cite{At}[Proposition 2.3].
\er

%%%%%%%%%%%%%%%%%%%%%%%%%%%%%

\subsection{Abstract Plethysm}

 Let $A$ be a pseudo-abelian (or Karoubian) $\bQ$-linear additive category which is also symmetric monoidal, with the tensor product $\otimes$ $\bQ$-linear additive in both variables. 
 \bd\label{seq} Let $S(A):=\prod_{n \geq 0} A_{\Sigma_n}$ (resp., $S_f(A):=\bigoplus_{n \geq 0} A_{\Sigma_n}$) be the {\it category of (finite) symmetric sequences} with values in $A$, endowed with the Cauchy product $\odot$ induced by ${\rm Ind}^{\Sigma_{n+m}}_{\Sigma_n \times \Sigma_m} : A_{\Sigma_n} \times A_{\Sigma_m} \to A_{\Sigma_{n+m}}$. When there is no danger of confusion, we will use the notation $S_{(f)}(A)$ to denote either of these categories.
\ed

The categories $S(A)$ and $S_f(A)$ are clearly pseudo-abelian $\bQ$-linear additive,  with the Cauchy product $\odot$ $\bQ$-linear additive in both variables. 
Moreover, these categories are also symmetric monoidal with respect to the Cauchy product $\odot$, with unit ${\bf 1}$ given by $1_A$ concentrated in degree zero (as in \cite{Fr1}[Sect.2.17]], \cite{Fr}[Sect.1.1.4], \cite{HSS}[Sect.2] or \cite{Ge}). Therefore, all results of Section \ref{simmon} are also applicable in the context of $S(A)$ and $S_f(A)$, with 
\be\label{se1} {\bar K}_0(S(A))=\prod_n {\bar K}^{\Sigma_n}_0(A)  \simeq {\bar K}_0(A) \otimes_{\bZ} \left( \prod_n Rep_{\bQ}(\Sigma_n) \right)= {\bar K}_0(A) \otimes_{\bZ}  {\bar K}_0 (S({\rm Vect}_{\bQ}))\ee and \be\label{se2}{\bar K}_0(S_f(A))=\bigoplus_n {\bar K}^{\Sigma_n}_0(A) \simeq {\bar K}_0(A) \otimes_{\bZ} \left( \bigoplus_n Rep_{\bQ}(\Sigma_n) \right)= {\bar K}_0(A) \otimes_{\bZ}  {\bar K}_0 (S_f({\rm Vect}_{\bQ})),\ee as ${\bar K}_0(A)$-algebras and commutative graded rings with the induced Cauchy products. 
In particular,  
both ${\bar K}_0(S(A))$ and ${\bar K}_0(S_f(A))$ are graded pre-lambda rings, and ${\bar K}_0(S(A))$ is the completion of ${\bar K}_0(S_f(A))$ with respect to the filtration induced by the grading.  The Adams operations on ${\bar K}_0(S_{(f)}(A))$ are identified via the above isomorphisms (\ref{se1}) and (\ref{se2}) with the tensor product of the corresponding Adams operations on the factor Grothendieck groups (which are ring homomorphisms by Lemma \ref{l1n}), i.e., for $\cV \in A$ and $V \in {\rm Vect}_{\bQ}(\Sigma_n)$, we have by 
$$(\cV \otimes V)^{\odot r} \simeq \cV^{\otimes r} \otimes V^{\odot r}$$
(as sequences concentrated in degree $r \cdot n$)
and by the multiplicativity of traces that
\be\label{Adcomp} \psi_r([\cV] \otimes [V])=\psi_r([\cV]) \otimes \psi_r([V]).\ee
In particular, $\psi_r$ is a graded ring homomorphism of degree $r$.

\br\label{r1n}
The above description of the Adams operations on ${\bar K}_0(A) \otimes_{\bZ} {\bar K}_0 (S_{(f)}({\rm Vect}_{\bQ}))$ would also follow from the fact that the ring isomorphisms (\ref{se1}) and (\ref{se2}) are in fact isomorphisms of pre-lambda rings, see \cite{Ge}[Theorem 4.8]. However, the latter fact is not needed in this paper.
\er

By using in (\ref{circ}) the category $S_{(f)}(A)$ instead of $A$, we get the {\it composition pairing} 
\be\label{circs} \circ: Rep_{\bQ}(\Sigma) \times {\bar K}_0(S_{(f)}A) \lra {\bar K}_0(S_{(f)}A)\ee corresponding to the action of polynomial functors on $S_{(f)}(A)$. By (\ref{se2}), this can be extended to a pairing
\be\label{circt} \circ:  {\bar K}_0(S_f A)  \times {\bar K}_0(S_{(f)}A) \lra {\bar K}_0(S_{(f)}A),\ee which is a ${\bar K}_0(A)$-algebra homomorphism in the first variable. By completion with respect to the degree filtration, we get the composition pairing 
\be\label{circu} \circ:  {\bar K}_0(S({\rm Vect}_{\bQ})) \times F_1 {\bar K}_0(S(A)) \lra {\bar K}_0(S(A))\ee corresponding to the action of analytic functors on $S(A)$. By (\ref{se1}), this can be further extended to a pairing
\be\label{circus} \circ:  {\bar K}_0(S(A))  \times F_1{\bar K}_0(S(A)) \lra {\bar K}_0(S(A)),\ee which is a ${\bar K}_0(A)$-algebra homomorphism in the first variable. Here, $F_1{\bar K}_0S(A)$ denotes the Grothendieck group of symmetric sequences with a zero entry in degree zero.

\medskip

We have from (\ref{cln}) the graded ring homomorphisms (with $p_i$ viewed as a variable in degree $i$)
$$cl:{\bar K}_0(S_f(A)) \lra {\bar K}_0(A) \otimes_{\bZ} \bQ[p_i, i \geq 1]$$
and, in particular, for $A={\rm Vect}_{\bQ}$,
$$cl: Rep_{\bQ}(\Sigma)={\bar K}_0(S_f({\rm Vect}_{\bQ})) \lra \bQ[p_i, i \geq 1].$$
By completion with respect to the grading, we also get
$$cl:{\bar K}_0(S(A)) \lra {\bar K}_0(A) \otimes_{\bZ} \bQ[[p_i, i \geq 1]],$$
mapping $F_1{\bar K}_0(S(A))$ to $F_1\bQ[[p_i, i \geq 1]]$, which  consists of power series with zero constant term.
In particular, for $A={\rm Vect}_{\bQ}$, we also get a graded ring homomorphism
$$cl: {\bar K}_0(S({\rm Vect}_{\bQ})) \lra \bQ[[p_i, i \geq 1]].$$

\bd\label{pl} {\it Abstract Plethysm} \\
Let $$f(p_1,p_2, \cdots, p_i,\cdots)=\sum_{n\geq 0}\sum_{\lambda \vdash n} f_{\lambda} \otimes p_{\lambda} \in {\bar K}_0(A) \otimes_{\bZ} \bQ[[p_i, i \geq 1]]$$
$$g(p_1,p_2, \cdots, p_j,\cdots)=\sum_{n\geq 0}\sum_{\mu \vdash n} g_{\mu} \otimes p_{\mu}\in {\bar K}_0(A) \otimes_{\bZ} \bQ[[p_i, i \geq 1]]$$ be two formal power series in the variables $p_i$ ($i \geq 1$), with coefficients $f_{\lambda} , g_{\mu} \in {\bar K}_0(A) \otimes_{\bZ} \bQ$. Assume that either $f$ is a polynomial in the $p_i$'s or $g$ has constant term $g_0=0$. The {\it plethysm} $f \circ g$ of $f$ and $g$ is defined by 
$$(f\circ g)(p_1,p_2, \cdots, p_i,\cdots):=$$
$$f\left(\psi_1g(p_1,p_2, \cdots, p_j,\cdots), \psi_2g(p_2,p_4, \cdots, p_{2j},\cdots),\cdots, \psi_ig(p_i,p_{2i} \cdots, p_{ij},\cdots), \cdots\right),$$
i.e., the variable $p_i$ in $f$ is substituted by the formal power series
$\psi_ig(p_i,p_{2i} \cdots, p_{ij},\cdots)$, with the power series $\psi_ig$ defined by the action of the Adams operation on the coefficients $g_{\mu}$ of $g$, i.e.,
$$\psi_ig(p_1,p_2, \cdots, p_j,\cdots):=\sum_{n\geq 0}\sum_{\mu \vdash n} \psi_i(g_{\mu}) \otimes p_{\mu}\in {\bar K}_0(A) \otimes_{\bZ} \bQ[[p_i, i \geq 1]].$$
Via the ring homomorphism $\bQ[[p_i, i \geq 1]] \lra {\bar K}_0(A) \otimes_{\bZ} \bQ[[p_i, i \geq 1]]$, one can also define the plethysm for $f \in \bQ[[p_i, i \geq 1]]$ and $g \in {\bar K}_0(A) \otimes_{\bZ} \bQ[[p_i, i \geq 1]]$ (with either $f$ a polynomial or $g$ with zero constant term).
\ed

\br In the special case $A={\rm Vect}_{\bQ}$, the corresponding  Adams operations act trivially on ${\bar K}_0({\rm Vect}_{\bQ})=\bZ$, so the above definition reduces to the classical notion of plethysm as in, e.g., \cite{J1,J2,LR}.
\er

The third abstract formula of this paper is given by the following generalization of a corresponding result of Joyal for the classical case $A={\rm Vect}_{\bQ}$ (see \cite{J1,J2}).
\bt\label{mt3} Let $\cV_{\cdot}=(\cV_n)_{n \geq 0}$ , $\cW_{\cdot}=(\cW_n)_{n \geq 0} \in S(A)$ be two symmetric sequences in $A$, with either $\cV_{\cdot}$ finite or $\cW_0=0_A$. Then we have
\be\label{ma3}
cl([\cV_{\cdot}] \circ [\cW_{\cdot}])= cl([\cV_{\cdot}] )\circ cl([\cW_{\cdot}]),
\ee
where on the left-hand side, $\circ$ is one the composition pairings (\ref{circt}) or (\ref{circus}), while on the right-hand side, $\circ$ denotes the abstract plethysm. A similar formula holds for $V_{\cdot}=(V_n)_{n \geq 0} \in S({\rm Vect}_{\bQ})$ and  $\cW_{\cdot}=(\cW_n)_{n \geq 0} \in S(A)$, with either $V_{\cdot}$ finite or $\cW_0=0_A$, where one the composition pairings (\ref{circs}) or (\ref{circu}) is used on the left-hand side.
\et

\begin{proof}
Since $cl$ and the composition product $\circ$ are compatible with the completion, it suffices to prove (\ref{ma3}) for bounded sequences.
Since $cl$ is a ${\bar K}_0(A)$-algebra homomorphism, and both the composition and the plethysm product are ${\bar K}_0(A)$-algebra homomorphisms in the first variable, we can further reduce the proof of  (\ref{ma3}) to the case $V_{\cdot}=(V_n)_{n \geq 0} \in S_f({\rm Vect}_{\bQ})$ and  $\cW_{\cdot}=(\cW_n)_{n \geq 0} \in S_f(A)$. Since $cl$ is a graded group homomorphism, and the composition and plethysm product are additive in the first variable, we can moreover assume that $V_{\cdot}$ is concentrated in one degree $n$, and it is given by a $\Sigma_n$-representation $V$.

With the above reductions, we have:
$$cl([V] \circ [\cW_{\cdot}])=cl(S_V([\cW_{\cdot}])).$$
By the definition of the composition pairing and (\ref{sub1}), 
\begin{equation*}
\begin{split} 
cl(S_V([\cW_{\cdot}]))&= cl\big( cl_n(V)(\psi_1([\cW_{\cdot}]), \cdots, \psi_i([\cW_{\cdot}]), \cdots) \big) \\
&= cl_n(V)(cl(\psi_1([\cW_{\cdot}])), \cdots, cl(\psi_i([\cW_{\cdot}])), \cdots),
\end{split}
\end{equation*}
where the last equality uses the fact that $cl$ is a ring homomorphism. It remains to show that for all $i \geq 1$, and $g:=cl([\cW_{\cdot}])(p_1,\cdots, p_j, \cdots) \in {\bar K}_0(A) \otimes_{\bZ} \bQ[p_j, j \geq 1]$, we have
\be\label{la}cl(\psi_i([\cW_{\cdot}]))(p_1,\cdots, p_j, \cdots) =\psi_ig(p_i,p_{2i} \cdots, p_{ij},\cdots).\ee
Since $cl$ and $\psi_i$ are additive, we can assume that 
$$[\cW_{\cdot}] =[\cW] \otimes [W_{\cdot}]  \in {\bar K}_0(A) \otimes_{\bZ} Rep_{\bQ}(\Sigma) \simeq  {\bar K}_0(S_f(A)),$$
with $g=cl([\cW_{\cdot}]) =[\cW] \otimes cl([W_{\cdot}]) \in {\bar K}_0(A) \otimes_{\bZ} \bQ[p_j, j \geq 1]$. 
Then we get by (\ref{Adcomp}) that
$$cl(\psi_i([\cW_{\cdot}]))=\psi_i([\cW]) \otimes cl(\psi_i([W_{\cdot}])).$$
Finally, by \cite{Mc}[Sect.I.8] or \cite{SS}[Thm.3.7, Thm.4.2], we have:
$$ cl(\psi_i([W_{\cdot}]))(p_1,\cdots, p_j, \cdots) =cl([W_{\cdot}])(p_i,p_{2i}, \cdots, p_{ij}, \cdots),$$ i.e., the Adams operation $\psi_i$ on $Rep_{\bQ}(\Sigma)$ corresponds under $cl$ to the $\bQ$-algebra homomorphism $p_j \mapsto p_{ij}$, $j \geq 1$. Note that \cite{SS} proves the above identity for Adams operations defined via traces (like in this paper). This completes the proof.
 \end{proof}
 
 \bex\label{pe} {\it Plethistic exponential isomorphism}\\
 Let $${\rm Exp}:=\sum_{n\geq 0} S_{triv_n}: \left(F_1{\bar K}_0(S(A)),\oplus\right) \overset{\simeq}{\lra} \left(1+ F_1{\bar K}_0(S(A)),\odot\right) $$
 be the plethistic exponential isomorphism, as in \cite{Ge}, with $S_{triv_n}$ the Schur functor associated to the trivial $\Sigma_n$-representation (which calculates the $n$-th symmetric power). We have
 $${\rm Exp}=triv \ \circ  \ : F_1{\bar K}_0(S(A)) \lra 1+ F_1{\bar K}_0(S(A)),$$
 where $triv:=(triv_n)_{n\geq 0}$ and with composition product $\circ$ given by (\ref{circu}). 
 Then
 $$cl(triv)=1 + \sum_{n \geq 1} h_n=\exp \left( \sum_{r \geq 1} \frac{p_r}{r} \right) \in \bQ[p_i, i \geq 1].$$ So, for any sequence $\cW_{\cdot}=(\cW_n)_{n \geq 0} \in S(A)$ with $W_0=0_A$, we have by Theorem \ref{mt3} that:
 $$cl({\rm Exp}([\cW_{\cdot}])=\exp \left( \sum_{r \geq 1} \frac{p_r}{r} \right) \circ cl([\cW_{\cdot}]).$$
 In particular, if $\cW_{\cdot}$ is a sequence concentrated in degree $1$ given by $\cW \in A$, we have that $S_{triv_n}(\cW)$ is concentrated in degree $n$, and it is given by $\cW^{\otimes n}$ with its natural $\Sigma_n$-action. Moreover, 
 $$cl([\cW_{\cdot}])=[\cW] \otimes p_1,$$
 so 
  $$cl({\rm Exp}([\cW_{\cdot}])=\exp \left( \sum_{r \geq 1} \frac{p_r}{r} \right) \circ ([\cW] \otimes p_1)=\exp \left(\sum_{r \geq 1} \psi_r([\cW]) \otimes \frac{p_r}{r} \right).$$
 Since ${\bar K}_0(S(A))=\prod_{n \geq 0} {\bar K}_0^{\Sigma_n}(A)$ is a graded ring, it can be identified with the (pre-lambda) subring  $$\sum_{n \geq 0} {\bar K}_0^{\Sigma_n}(A) \cdot t^n\subset {\bar K}_0(S_f(A))[[t]],$$  
 and similarly for ${\bar K}_0(A) \otimes_{\bZ} \bQ[[p_i, i \geq 1]]$.
 Then the above formula recovers (\ref{mainabsii}), i.e., 
$$cl(\sum_{n\geq 0} [\cW^{\otimes n}] \cdot t^n)=\exp \left(\sum_{r \geq 1} \psi_r([\cW]) \otimes p_r \cdot \frac{t^r}{r} \right),$$ as $p_r$ is a variable of degree $r$ (see also \cite{Ge}[Prop.5.4]).
 \eex
 
The Grothendieck group ${\bar K}_0(S_{(f)}(A))$ has more structure than needed in this paper (except that Remark \ref{plinv} is used in the Example \ref{hilbord2} from the introduction). This will be explained in the following remarks. 

 \br\label{plinv} 
 The composition products (\ref{circs})--(\ref{circus}) are associative (see \cite{Ge,GP}), which implies that the ring ${\bar K}_0(S_{(f)}(A))$ and the degree-zero subring ${\bar K}_0(A)$ are lambda rings (see \cite{Ge}[Thm.4.7], and compare also with \cite{He}). In particular, the Adams operations $\psi_i:{\bar K}_0(A) \to {\bar K}_0(A)$ satisfy (in addition to Lemma \ref{l1n}) the relation $\psi_i \circ \psi_j=\psi_{ij}$, for any positive integers $i,j$. 

Similar results apply to the usual Grothendieck groups
$K_0(S_{(f)}(A))$ and the degree-zero subring $K_0(A)$ in case $A$ (and therefore also $S_{(f)}(A)$) is an {\it Abelian category} with the tensor product {\it exac}t in both variables. Then the natural maps ${\bar K}_0(S_{(f)}(A))\to K_0(S_{(f)}(A))$ and
${\bar K}_0(A)\to K_0(A)$ are (pre-) lambda ring homomorphisms. Moreover, the first one commutes with the corresponding
composition product $\circ$ and the
plethistic exponential isomorphism ${\rm Exp}$ (see \cite{Ge, GP}).

Associativity is also used in \cite{Ge,H} to reduce the calculation of invariants of configuration spaces $\bF(X,n)$ and  their Fulton-MacPherson compactifications $\bF\bM(X,n)$ (for $X$ smooth) to the calculation of the corresponding invariants for products, namely:
 \be\label{pl1} \left([H^*_c(X^n)]\right)_{n \geq 0}=\left([H^*_c(\bF(X,n))]\right)_{n \geq 0} 
 \circ (0, triv_{1}, \cdots, triv_j, \cdots), \ee 
 resp., 
  \be\label{pl2} \left([H^*_c(\bF\bM(X,n))]\right)_{n \geq 0}=\left([H^*_c(\bF(X,n))]\right)_{n \geq 0} 
 \circ ([H^*({\rm P}_n)])_{n \geq 0},\ee  
 for a suitable sequence $({\rm P}_n)_{n \geq 0}$ of $\Sigma_n$-equivariant smooth projective varieties, as in \cite{Ge,H}
(and with ${\rm P}_n$  depending only on $\dim(X)$, but not on $X$, and  ${\rm P}_0:=\emptyset$). Here, $[H^*_c]:=\sum_{i} (-1)^i [H^i_c]\in K_0$ in the usual Grothendieck group $K_0$ of the corresponding Abelian category.
Moreover, it is important to work with the cohomology with compact support $H^*_c$ and the Euler-characteristic type class 
$[H^*_c]$, since ``additivity'' with respect to decompositions is used in these calculations
(i.e., $[H^*_c]$ is a ``Serre functor'' in the sense of \cite{Ge}[Def.5.2]).
 
 Note that  $$(0, triv_{1}, \cdots, triv_j, \cdots)={\rm Exp}(triv_1)-1 \in {\bar K}_0(S({\rm Vect}_{\bQ}))
= K_0(S({\rm Vect}_{\bQ})),$$ 
and the plethistic logarithm isomorphism (inverse to ${\rm Exp}$) is given by
 $${\rm Log(1+\cdot}):= \left( {\rm Exp}(triv_1)-1 \right)^{-\circ}  \circ \cdot : F_1K_0(S(A)) \overset{\simeq}{\lra} F_1K_0(S(A)),$$
(see \cite{Ge}[Prop.2.2], \cite{GP}[Sec.5]). Here, $(-)^{-\circ}$ denotes the inverse with respect to composition. Therefore, by using the associativity of the composition, (\ref{pl1}) and (\ref{pl2}) can be reformulated as (for a ``motivic'' version of the first formula  see also \cite{GP}[Thm. 3.2]):
 \be\label{pl1b}
 \left([H^*_c(\bF(X,n))]\right)_{n \geq 0}={\rm Exp} \Big( [H^*_c(X)] \otimes {\rm Log}(1+triv_1) \Big)=:(1+triv_1)^{[H^*_c(X)]}
 \ee
 and 
   \be\label{pl2b} 
   \begin{split} \left([H^*_c(\bF\bM(X,n))]\right)_{n \geq 0} &={\rm Exp} \left( [H^*_c(X)] \otimes {\rm Log}\Big(1+\sum_{i \geq 1} [H^*_c({\rm P}_i)]\Big) \right).\\
   &=:\Big(1+\sum_{i \geq 1} [H^*_c({\rm P}_i)]\Big)^{[H^*_c(X)]} \in K_0(S(A)).
   \end{split}
 \ee
Here $\otimes$ can be  the composition $\circ$, with   $[H^*_c(X)]$ viewed as symmetric sequence concentrated in degree one, or equivalently, the Cauchy product $\odot$, with   $[H^*_c(X)]$ viewed as symmetric sequence concentrated in degree zero.
\er

\br\label{Had} 
\noindent{(1)} The pairing (\ref{pair}) can be degree-wise extended to a pairing (also called the {\it Hadamard product}):
 $$\ast: {\bar K}_0(S_{(f)}({\rm Vect}_{\bQ})) \times {\bar K}_0(S_{(f)}(A)) {\lra} {\bar K}_0(S_{(f)}(A));\;
[V_{\cdot}] \ast [\cW_{\cdot}]:=[(V_{n} \otimes  \cW_{n})_{n \geq 0}]$$ such that  
 \be\label{ma5}
cl([V_{\cdot}] \ast [\cW_{\cdot}])= cl([V_{\cdot}] )\ast cl([\cW_{\cdot}]),
\ee
for any symmetric sequences $V_{\cdot}=(V_n)_{n \geq 0} \in S_{(f)}({\rm Vect}_{\bQ})$ and  $\cW_{\cdot}=(\cW_n)_{n \geq 0} \in S_{(f)}(A)$. The product on the right-hand side of (\ref{ma5}) is  the {\it internal product} of formal power series in the $p_i$ ($i \geq 1$), which is defined by the rule $p_{\lambda} \ast p_{\mu} = \delta_{\lambda\mu} z_{\lambda} p_{\lambda}$, as in \cite{Mc}[(I.7.12)]. If, for example, we take $\cW_n:=\cV^{\otimes n}$, for some $\cV \in A$, then (\ref{ma5}) generalizes  formula (\ref{new3a}) of Theorem \ref{abstwist}.

A  formula similar to (\ref{ma5}) holds for the Hadamard product 
 $$\ast: {\bar K}_0(S_{(f)}(A))  \times {\bar K}_0(S_{(f)}(A)) {\lra} {\bar K}_0(S_{(f)}(A)));\;
[\cV_{\cdot}] \ast [\cW_{\cdot}]:=[(\cV_{n} \otimes  \cW_{n})_{n \geq 0}]$$
 induced degree-wise from the symmetric monoidal structure.

\noindent{(2)} The algebra ${\bar K}_0(S(A))$ is a {\it complete algebra with composition} operation $\circ$ over ${\bar K}_0(A)$, in the sense of \cite{GP}, see Theorem 5.1 in loc.cit.
 \er

%%%%%%%%%%%%%%%%%%%%%%%%%%%%%
%%%%%%%%%%%%%%%%%%%%%%%%%%%%%

\section{From abstract to concrete identities. Examples and Applications}
In this section, we show how to derive Theorem \ref{mth} and formula (\ref{new1}) from the Introduction as consequences of the abstract general results of the previous section. Furthermore, we illustrate various special cases of Theorems \ref{mth} and \ref{tw} in Subsection \ref{spc}. Finally, in Subsection \ref{feq} we present further applications of the abstract formulae to a suitable equivariant context.

\subsection{From abstract to concrete identities.}\label{pmt}
Let us now explain how to derive our Theorem \ref{mth}, as well as formula (\ref{new1}) from the Introduction from the abstract generating series formula (\ref{mainabsii}) of Section \ref{simmon}. We start with the proof of formula (\ref{mf1h}) in the mixed Hodge context. 

For an additive tensor category $(Ab,\otimes)$,  let $Gr^{-}(Ab)$  denote the additive tensor category of {bounded graded}
objects in $Ab$, i.e., functors $G: \bZ\to Ab$, with $G_n:=G(n)=0$ except for finitely many $n\in \bZ$. Here, $$(G\otimes G')_n:= \bigoplus_{i+j=n} G_i\otimes G_j \:,$$
with the Koszul symmetry isomorphism (indicated by the $-$ sign in $Gr^{-}$): $$(-1)^{i\cdot j} s(G_i,G_j): G_i\otimes G_j \simeq  G_j\otimes G_i \:.$$
If $(Ab,\otimes)$ is a $\bQ$-linear Karoubian (or abelian)  symmetric monoidal category, then the same is true for $Gr^{-}(Ab)$. This applies for example to the category $\ms$ of mixed Hodge structures. Note that in the K\"unneth formula (\ref{Kue}), we have to view $H^*_{({c})}(X,\cM)$ as an element in the $Gr^{-}(\ms)$ with tensor product $\otimes$ defined via the above Koszul rule.

Let $Gr_F^*Gr_*^W: \ms\to  Gr^2({\rm Vect}_{\bC})$ be the functor of taking the associated bigraded finite dimensional $\bC$-vector space:
$$V \mapsto \bigoplus_{p,q} \; Gr_F^pGr^W_{p+q}(V \otimes_{\bQ}\bC) \in  Gr^2({\rm Vect}_{\bC}) \:.$$
This is an exact tensor functor of such abelian tensor categories, if we use the induced symmetry isomorphism without any sign changes for the abelian category $Gr^2({\rm Vect}_{\bC})$ of bigraded finite dimensional complex vector spaces.
The transformation  $Gr_F^*Gr_*^W$ is compatible with the K\"unneth isomorphism (\ref{Kue}). Similarly, $Gr_F^*Gr_*^W$ is compatible with the abstract pairing (\ref{pair}), as well as taking invariant subobjects. Moreover, for $A=Gr^{-}(\ms)$, the abstract pairing gets identified with the tensor product on $A$, as used in (\ref{tKue}), after regarding a rational representation as a pure Hodge structure of type $(0,0)$ placed in degree zero.

Recall next that the ring homomorphism $$h: {K}_0(Gr^{-}(Gr^2({\rm Vect}_{\bC})))\to \bZ[y^{\pm 1},x^{\pm 1},z^{\pm 1}] $$
given by $$\left[\bigoplus (V^{p,q})^k\right] \mapsto \sum_{p,q,k} \;dim((V^{p,q})^k) \cdot y^px^q(-z)^k \:,$$
with $k$ the degree with respect to the grading in $Gr^{-}$ is an isomorphism of pre-lambda rings, see \cite{MS}[Prop.2.4]. The pre-lambda structure on ${K}_0(Gr^{-}(Gr^2({\rm Vect}_{\bC})))$ is defined as in (\ref{prel}), whereas the pre-lambda structure on the Laurent polynomial ring $\bZ[y^{\pm 1},x^{\pm 1},z^{\pm 1}]$ corresponds to the Adams operations 
$$\psi_r(p(y,x,z))=p(y^r,x^r,z^r).$$
The sign choice of numbering by $(-z)^k$ in the
definition of $h$ is needed for the compatibility with these pre-lambda structures.

Finally, we have an equality 
$$(h \otimes id) \circ tr_n=tr_{\Sigma_n} :K_0^{\Sigma_n}(Gr^{-}Gr^2({\rm Vect}_{\bC})) \to  \bZ[y^{\pm 1},x^{\pm 1},z^{\pm 1}]
 \otimes_{\bZ} C(\Sigma_n),$$
as can be easily checked on generators given by a $\Sigma_n$-representation placed in a single multi-degree.

Formula (\ref{mf1h}) follows now by applying the ring homomorphism 
$$(h\circ Gr_F^*Gr_*^W)  \otimes id : K_0(Gr^{-}(\ms)) \otimes_{\bZ} \bQ[p_i, i \geq 1] \to \bZ[y^{\pm 1}, x^{\pm 1},z^{\pm 1}] \otimes_{\bZ} \bQ[p_i, i \geq 1]$$
to formula (\ref{mainabsii}) of Theorem \ref{mta}, with $\cV:=H^*_{({c})}(X,\cM) \in A:=Gr^{-}(\ms)$.
Similarly, formula (\ref{new1}) follows by applying this ring homomorphism to the identity (\ref{npairing}).

\medskip

For the proof of the generating series (\ref{mf1}) and the multiplicativity (\ref{new1}) for the Poincar\'e-type polynomials, we consider similarly the isomorphism of pre-lambda rings
$$P: K_0(Gr^{-}({\rm Vect}_{\bK}))\to \bZ[z^{\pm 1}]$$
defined by taking the dimension counting Laurent polynomial
$$\left[\bigoplus V^k\right] \mapsto \sum_{k} \;dim(V^k) \cdot (-z)^k \:,$$
with $k$ the degree with respect to the grading in $Gr^{-}$.
Here, ${\rm Vect}_{\bK}$ is the abelian tensor category of finite dimensional $\bK$- vector spaces (with $\bK$ a field of characteristic zero), and the Adams operation on $\bZ[z^{\pm 1}]$ is given by $\psi_r(p(z))=p(z^r).$
Similarly, we have an equality 
$$(P \otimes id) \circ tr_n=tr_{\Sigma_n} :K_0^{\Sigma_n}(Gr^{-}({\rm Vect}_{\bK})) \to  \bZ[z^{\pm 1}]\otimes_{\bZ} C(\Sigma_n)  .$$
Then formula (\ref{mf1}) follows by applying the ring homomorphism 
$$P  \otimes id : K_0(Gr^{-}({\rm Vect}_{\bK})) \otimes_{\bZ} \bQ[p_i, i \geq 1] \to \bZ[z^{\pm 1}] \otimes_{\bZ} \bQ[p_i, i \geq 1]$$
to formula (\ref{mainabsi}), with $\cV:=H^*_{({c})}(X,\cM) \in A:=Gr^{-}({\rm Vect}_{\bK})$. Similarly, formula (\ref{new1}) follows by applying this ring homomorphism to the identity (\ref{npairing}).

%%%%%%%%%%%%%%%%%%%%%%%%%%%%%%%%
%%%%%%%%%%%%%%%%%%%%%%%%%%%%%%%%

\subsection{Examples}\label{spc}
The results of Theorems \ref{mth} and \ref{tw} can be specialized in several different ways, e.g., 
\begin{itemize}
\item[(i)] for specific values of the parameter $z$ (and, resp., $x$, $y$, $z$ in the mixed Hodge context), e.g., the specialization $z=1$ yields Euler-characteristic type formulae; \item[(ii)] for special choices of the coefficients $\cM$ (e.g., intersection cohomology complexes); \item[(iii)] for special values of the Frobenius parameters $p_r$ (e.g., related to symmetric and alternating powers of coefficients).
\item[(iv)] for special choices of the representation $V \in Rep_{\bQ}(\Sigma_n)$, e.g., for $V={\rm Ind}_K^{\Sigma_n}(triv_K)$, the representation induced from the trivial representation of a subgroup $K$ of $\Sigma_n$.
\end{itemize}

For the convenience of the reader, we next explain some of these cases in detail.

\medskip

\noindent $(i)$ \ By letting $z=1$ in (\ref{mf1}) we obtain a generating series identity for the characters of the virtual cohomology representations $$[H_{({c})}^*(X^n,\cM^{\boxtimes n})]:=\sum_k (-1)^k [H_{({c})}^k(X^n,\cM^{\boxtimes n})] \in Rep_{\bK}(\Sigma_n),$$ with $P_{({c})}$ on the right-hand side of (\ref{mf1}) being replaced by the corresponding (compactly supported) {\it Euler characteristic} $$\chi_{({c})}(X,\cM):=\sum_{k} (-1)^k \cdot b_{({c})}^k(X,\cM) \in \bZ.$$ Similarly, 
by letting $z=1$ in (\ref{mf1h}), we get a generating series formula for the characters of graded parts (with respect to both filtrations) of the virtual cohomology representations 
$$\sum_{k,p,q} (-1)^k \cdot [Gr^p_FGr_{p+q}^WH_{({c})}^k(X^n,\cM^{\boxtimes n})] y^px^q \in Rep_{\bC}(\Sigma_n)[y^{\pm 1}, x^{\pm 1}],$$ where $h_{({c})}$ on the right-hand side of (\ref{mf1h}) gets replaced by its specialization to the 
{\it E-polynomia}l $E_{({c})}$. In this case, we recast Getzler's generating series formula for the $E$-polynomial \cite{Ge}[Prop.5.4].
Finally, by letting $x=z=1$ in (\ref{mf1h}), we get a generating series formula for the 
characters of graded parts (with respect to the Hodge filtration) of the virtual cohomology representations 
$$\sum_{k,p} (-1)^k \cdot [Gr^p_FH_{({c})}^k(X^n,\cM^{\boxtimes n})] y^p \in Rep_{\bC}(\Sigma_n)[y^{\pm 1}],$$ where $h_{({c})}$ on the right-hand side of (\ref{mf1h}) gets replaced by its specialization to the Hodge polynomial (or {\it Hirzebruch characteristic}) $\chi^{(c)}_{-y}$.
\br
If $X$ is a {\it complex projective} variety, some of these special cases of Euler characteristic-type generating series have been derived in our recent work \cite{MS15}[Eqn.(9),(10)] by taking degrees of suitable equivariant characteristic class formulae. Note that in the mixed Hodge context, these characteristic class formulae only take into account the Hodge filtration, so the $E$-polynomial version discussed above, as well as Theorem \ref{mth} cannot be deduced as degree formulae. Moreover, if $X$ is a complex projective manifold, the specialization $\chi_{y}$ mentioned above becomes the classical {\it Hirzebruch $\chi_y$-genus}. 
\er

%\medskip

\noindent  $(ii)$ \ Let us now specialize Theorem \ref{mth} to the important concrete example of coefficients (other than trivial coefficients) given by the intersection cohomology complex, e.g., in the Hodge context.
Let $X$ be a pure-dimensional complex quasi-projective variety, and let $\cM$ be the (shifted) {\it intersection cohomology Hodge module} $${IC'}_X^H:=IC^H_X[-\dim(X)] \in D^b\mh(X),$$ with underlying constructible sheaf complex
 $IC'_X:=IC_X[-\dim(X)] \in D_c^b(X)$.
The (compactly supported) cohomology groups $H_{({c})}^*(X,{IC'}_X^H)$ endow the (compactly supported) intersection cohomology groups of $X$, that is, 
$$IH^*_{({c})}(X):=H_{({c})}^*(X,IC'_X),$$
with mixed Hodge structures.
Thus, as a special case of  (\ref{mf1h}) we get a generating series formula for the characters of graded parts (with respect to both filtrations) of {\it intersection cohomology representations} of cartesian products of $X$, namely:
\be\label{mc1ih}
\sum_{n \geq 0} tr_{\Sigma_n}(IH^*_{({c})}(X^n)) \cdot t^n  
 =\exp \left( \sum_{r \geq 1}  p_r \cdot {h_{({c})}}(X,{IC'}_X^H)(y^r,x^r,z^r) \cdot \frac{t^r}{r}\right).
\ee
By letting $y=x=1$ in (\ref{mc1ih}), we obtain a generating series for the corresponding Poincar\'e-type polynomials of characters of intersection cohomology representations of cartesian products of $X$.\\

\noindent $(iii)$ \ For suitable values of the Frobenius parameters $p_r$ in Theorem \ref{mth}, formulae (\ref{mf1}) and  (\ref{mf1h}) also generalize several generating series identities from \cite{MS} for the Betti numbers (respectively, mixed Hodge numbers) of {\it symmetric} and resp. {\it alternating powers} $\cM^{(n)}$, resp., $\cM^{\{n\}}$ of the coefficient  $\cM$; see Section \ref{psf} for a definition of these objects, which are induced coefficients on the symmetric product $X^{(n)}:=X^n/{\Sigma_n}$. 

More precisely, by making $p_r=1$ for all $r$, the effect is to take the $\Sigma_n$-invariant part in the K\"unneth formula, i.e., to compute the Betti (or mixed Hodge) numbers of
\be\label{kea} H_{({c})}^*(X^n,\cM^{\boxtimes n})^{\Sigma_n} \cong H^*_{({c})}(X^{(n)},\cM^{(n)}).\ee
E.g., we recover from (\ref{mf1}) the generating series for the Poincar\'e polynomials and Betti numbers of symmetric powers $\cM^{(n)}$ of $\cM$, namely:
\be\label{i8}
 \sum_{n \geq 0}  P_{({c})}(X^{(n)},\cM^{(n)})(z) \cdot t^n  
=\exp \left( \sum_{r \geq 1}  P_{({c})}(X,\cM)(z^r) \cdot \frac{t^r}{r} \right),
\ee 
and similarly for the mixed Hodge polynomials in the mixed Hodge context. For more applications and special cases of formula (\ref{i8}) and its Hodge version, the reader is advised to consult our previous work \cite{MS}. Let us only recall here that for the trivial coefficients $\cM = \bK_X, \ \cO_X, \ \bQ^H_X$, we have that $\cM^{(n)}=\bK_{X^{(n)}}, \ \cO_{X^{(n)}},, \ \bQ^H_{X^{(n)}},$ and similarly for the (shifted) intersection cohomology Hodge modules, i.e.,
$({IC'}_X^H)^{(n)}={IC'}^H_{X^{(n)}}.$
 
By making $p_r=(-1)^{r-1}$ for all $r$,  the effect is to take the $\Sigma_n$-anti-invariant part in the K\"unneth formula, i.e., to compute the Betti (or mixed Hodge) numbers of
\be\label{keb} H_{({c})}^*(X^n,\cM^{\boxtimes n})^{sign-\Sigma_n} \cong H^*_{({c})}(X^{(n)},\cM^{\{n\}}).\ee
So we obtain a generating series formula for the Betti numbers $b_{({c})}^{k}(X^{(n)},\cM^{\{n\}})$ and, respectively, mixed  Hodge numbers $h_{({c})}^{p,q,k}(X^{(n)},\cM^{\{n\}})$ if $\cM \in D^b\mh(X)$, of the alternating powers of $\cM$, i.e., 
  \be\label{cf}
  \sum_{n \geq 0}  P_{c}(X^{(n)},\cM^{\{n\}})(z) \cdot t^n  
=\exp \left( \sum_{r \geq 1} - P_{c}(X,\cM)(z^r) \cdot \frac{(-t)^r}{r} \right),
 \ee
and  similarly for the mixed Hodge polynomials in the mixed Hodge context.
For concrete examples and special cases of formula \ref{cf}, see \cite{MS}.
 
 \medskip
 
The specialization $p_1 \mapsto 1$ and $p_r \mapsto 0$ if $r \geq 2$ corresponds (by (\ref{defF})) to forgetting the $\Sigma_n$-action, up to the  Frobenius-type factor $\frac{1}{n!}$.
So, as a consequence of Theorem \ref{mth}, we also get the following:
 \bc\label{c13} Under the notations and assumptions of Theorem \ref{mth}, 
  the following generating series holds in $\bQ[z^{\pm 1}][[t]]$:
 \be\label{K1}
 \sum_{n \geq 0} P_{({c})}(X^n,\cM^{\boxtimes n})(z) \cdot \frac{t^n}{n!}  
 =\exp \left( {P_{({c})}}(X,\cM)(z) \cdot t \right),
 \ee
 and similarly for the mixed Hodge polynomials in the mixed Hodge context.\ec
 
Formula (\ref{K1}) and its Hodge theoretic version can also be obtained directly from the K\"unneth isomorphisms (\ref{Kut}) or (\ref{Kue}), depending on the context.\\

\noindent $(iv)$ \ 
Interesting new specializations of our results also arise for different choices of the representation $V$. For example, by choosing $V={\rm Ind}_K^{\Sigma_n}(triv_K)$, the representation induced from the trivial representation of a subgroup $K$ of $\Sigma_n$, and for $\cM=\bK_X \in D^b_c(X;\bK)$ the constant sheaf, formula (\ref{new2a})  specializes for $p_r=1$ (for all $r$) to Macdonald's Poincar\'e polynomial formula \cite{Mac}[p.567] for the quotient $X^n/K$, i.e.,
\be\label{pquot}
P_{({c})}(X^n/K)(z)= \sum_{{\lambda=(k_1,k_2, \cdots) \vdash n} } \frac{1}{z_{\lambda}} \chi_{\lambda}({\rm Ind}_K^{\Sigma_n}(triv_K))  \cdot \prod_{r \geq 1} \left( P_{({c})}(X)(z^r) \right)^{k_r} \in \bK[z].
\ee
A similar formula also holds for the mixed Hodge polynomial $h_{({c})}(X^n/K)(y,x,z)$ in the mixed Hodge context (see also (\ref{hquot})).
If $X$ is projective, similar identities hold for the Poincar\'e polynomial of the coherent structure sheaf $\cO_X$.

\medskip

Similarly, for $V$ a $\Sigma_n$-representation, (\ref{new2a}) specializes for $p_r=1$ (for all $r$) to a formula for the Poincar\'e (resp. mixed Hodge polynomials in the mixed Hodge context) of the corresponding {\it Schur-type object} $S_{V}(\cM)$ on $X^{(n)}$  (see Section \ref{psftwist} for a definition), namely:
\be\label{pschur}
P_{({c})}(X^{(n)},S_{V}(\cM))(z)= \sum_{{\lambda=(k_1,k_2, \cdots) \vdash n} } \frac{1}{z_{\lambda}} \chi_{\lambda}(V)  \cdot \prod_{r \geq 1} \left( P_{({c})}(H^*(X;\cM)(z^r) \right)^{k_r},
\ee
and similarly for the mixed Hodge context.
Note that at the cohomology level, we have the isomorphisms:
\be\label{sc}
H^*_{({c})}(X^{(n)},S_V(\cM)) \cong \left( V \otimes H^*_{({c})}(X^n,\cM^{\boxtimes}) \right)^{\Sigma_n},
\ee
and similarly for the graded pieces with respect to the Hodge and weight filtrations in the Hodge context.
These Schur-type objects $S_V(\cM)$ generalize the symmetric and alternating powers of $\cM$, which correspond to the trivial and, resp., sign representation. 

As a more concrete example, for $V=V_{\mu} \simeq V^*_{\mu}$ the (self-dual)  irreducible representation of $\Sigma_n$ corresponding to a partition $\mu$ of $n$, and  for $X$ a pure dimensional complex quasi-projective variety with $\cM=IC'_X \in D_c^b(X;\bQ)$, 
the corresponding Schur-type object $S_{V_\mu}(IC'_X)$ is given by the (shifted) twisted intersection cohomology complex $IC'_{X^{(n)}}(V_{\mu}) \in D^b_c(X;\bQ)$, with twisted coefficients corresponding to the local system on the configuration space $X^{\{n\}}\subset X^{(n)}$ of unordered $n$-tuples of distinct points in $X$, induced from $V_{\mu}$ by the group homomorphism $\pi_1(X^{\{n\}}) \to \Sigma_n$ (compare \cite{MS}[p.293] and \cite{MR}[Prop.3.5]). So, formula (\ref{pschur}) reduces in this case to the calculation of the Poincar\'e polynomials of twisted intersection cohomology $$IH_{({c})}^*(X^{(n)},V_{\mu}):=H_{({c})}^*(X^{(n)}, IC'_{X^{(n)}}(V_{\mu})),$$ namely, 
\be
P_{({c})}(X^{(n)},IC'_{X^{(n)}}(V_{\mu}))(z)= \sum_{{\lambda=(k_1,k_2, \cdots) \vdash n} } \frac{1}{z_{\lambda}} \chi_{\lambda}(V_{\mu})  \cdot \prod_{r \geq 1} \left( P_{({c})}(X, IC'_X)(z^r) \right)^{k_r}.
\ee
A similar formula holds for the mixed Hodge polynomials in the mixed Hodge context, i.e., for $\cM={IC'}_X^H \in D^b\mh(X)$ the (shifted) intersection cohomology Hodge module.
A special case of the latter, for the $\chi_{-y}$-polynomial, has been recently obtained by the authors in \cite{MS15}[Eqn.(23)], by taking degrees of a certain characteristic class identity.

%%%%%%%%%%%%%%%%%%%%%%%%%%%%%%%%%%%%%%
%%%%%%%%%%%%%%%%%%%%%%%%%%%%%%%%%%%%%%

\subsection{Further applications: equivariant context}\label{feq}
In this section, we indicate further applications of the abstract setup of the previous sections to suitable equivariant versions of (characters of) Poincar\'e and mixed Hodge polynomials of {\it equivariant} coefficients. 
More precisely, let $X$ and $H^*$ be as in Section \ref{gs} of the Introduction, with the $\Sigma_n$-equivariant formula (\ref{Kut}) being functorial. Moreover, we consider the following equivariant situations:
\begin{enumerate}
\item[(A)] $G$ is a fixed finite group acting (algebraically) on $X$;
\item[(B)] $g$ is a finite order (algebraic) automorphism acting on $X$;
\item[(C)] $g:X \to X$ is  an (algebraic) endomorphism, which moreover is required to be proper if compactly supported cohomology is used.
\end{enumerate} 
To unify the treatment of all these cases, note that any of the above situations can be viewed in the context of a (semi-)group
 action of $G$, with $G:=\bZ$ for (B) and $g=1 \in \bZ$ acting with finite order, and, resp., $G:=\bN_0$ for (C).

Our results and techniques also apply to  suitable $G$-equivariant coefficients $\cM$ on such a $G$-space $X$, 
provided that the finite dimensional $\mathbb{K}$-vector space $H^*_{({c})}(X,\cM)$ has an induced $G$-action, and a corresponding $G \times \Sigma_n$-equivariant K\"unneth formula holds:
\be\label{Kueq} H^*_{({c})}(X^n,\cM^{\boxtimes n}) \simeq H^*_{({c})}(X,\cM)^{\otimes{n}},\ee
where $\cM^{\boxtimes n}$ is the $n$-th self-external product of $\cM$ with its induced $\Sigma_n$-action commuting with the induced diagonal $G$-action.
Examples of such $G$-equivariant coefficients include objects $\cM\in A(X)$, with $A(X)$ any of the categories $D^b_{f}(X;\bK)$, $D^b_{c}(X;\bK)$, $D^b_{\it coh}(X)$ and $D^b\mh(X)$, and with the following types of $G$-action: 
\begin{enumerate}
\item[(A)] $G$ is a fixed finite group acting (algebraically) on $X$, with $\cM \in A_G(X)$ a $G$-equivariant object in $A(X)$ (as in \cite{MS}[Appendix]).
\item[(B)] $g$ is a finite order (algebraic) automorphism acting on $X$, with $\cM \in A_{\langle g \rangle}(X)$ a $\langle g \rangle$-equivariant object (in particular, $\cM \in A(X)$ is endowed with an isomorphism $\Psi_g:\cM \to g_*\cM$ in $A(X)$). Here the order of the cyclic group $\langle g \rangle$ can depend on $\cM$ (i.e., this order could exceed that of the action on $X$).
\item[(C)] $g:X \to X$ is  an (algebraic) endomorphism (which is required to be proper if compactly supported cohomology is used), together with a morphism $\Psi_g:\cM \to g_{*(!)}\cM$ in $A(X)$. Here, we use the derived pushforward $g_*$ (resp. $g_!$) when considering  (compactly supported) cohomology $H^*_{({c})}(X,\cM)$ with the endomorphism  induced from $\Psi_g$. Note that $g_!=g_*$ if $g$ is proper, e.g., an automorphism.
\end{enumerate} 
Concrete examples of such $G$-equivariant coefficients on a $G$-space $X$ include the ``trivial'' coefficients (i.e., the constant sheaf $\bK_X$,  the structure sheaf  $\cO_X$, and the constant Hodge sheaf $\bQ^H_X$). Here $\Psi_g$ is induced by the adjunction map $id\lra g_*g^*$
corresponding to the usual (derived) pullback in cohomology.
Similarly, in cases (A) and (B) one can use the intersection cohomology (Hodge) sheaf $IC'_X$ if $X$ is pure dimensional.

\medskip

For a $G$-equivariant object $\cM \in A(X)$ as above, the external products $\cM^{\boxtimes n}\in A(X^n)$  and their cohomology $H^*_{({c})}(X^n, \cM^{\boxtimes n})$
get an induced diagonal $G$-action commuting with the action of the symmetric group $\Sigma_n$ as before, so that for $V$ a $\Sigma_n$-representation (with trivial $G$-action), the (twisted) cohomology $H^*_{({c})}(X^n, V \otimes \cM^{\boxtimes n})$ also has an induced action of $G \times \Sigma_n$. In the mixed Hodge context, the induced action on $H^*_{({c})}(X^n, V \otimes \cM^{\boxtimes n})$ is an action in the category of rational mixed Hodge structures (since it is induced from a $G$-equivariant object in $D^b\mh(X)$). In particular, the $G$-action is compatible with the Hodge and, resp., weight filtration.

\medskip

All our concrete results from Sections \ref{gs} and \ref{twist} can be now formulated in the above equivariant context, once we redefine the Poincar\'e and resp. mixed Hodge polynomials, and the corresponding characters $tr_{\Sigma_n}$. Let
 $Rep_{\bK}(G):={\bar K}_0(A_G)$ denote the Grothendieck ring of the following $\bQ$-linear Karoubian (even abelian) tensor category $A_G$, corresponding to each of our situations above (with $\bK$ a field of characteristic zero):
 \begin{enumerate}
\item[(A)] ${\rm Vect}_{\bK}(G)$, the category of finite-dimensional $G$-representations.
\item[(B)] ${\rm Vect}_{\bK}(G)$, the category of finite-dimensional $G$-representations, with $g=1 \in G:=\bZ$ acting with finite order.
\item[(C)] ${\rm End}_{\bK}$, the category of endomorphisms of finite-dimensional $\bK$-vector spaces.
 \end{enumerate} 
  The tensor structure on $A_G$ is induced from the tensor product of the underlying $\bK$-vector spaces with induced diagonal $G$-action. Then a $G\times \Sigma_n$-action on a finitely dimensional vector space $V$ is the same a $\Sigma_n$-action on $V$ regarded as an object in $A_G$. By the Schur functor decomposition (\ref{Sch}) applied to $A_G$, we get the isomorphism $${\bar K}_0^{\Sigma_n}(A_G) \simeq Rep_{\bQ}(\Sigma_n) \otimes_{\bZ} {\bar K}_0(A_G)=Rep_{\bQ}(\Sigma_n) \otimes_{\bZ} Rep_{\bK}(G).$$ Then we define $G$-equivariant characters $$tr^G_{\Sigma_n}:{\bar K}_0^{\Sigma_n}(A_G) \to C(\Sigma_n) \otimes_{\bZ} Rep_{\bK}(G)$$ by taking the $\Sigma_n$-character in the first tensor factor of the above decomposition. 
 
We can now define {\it $Rep_{\bK}(G)$-valued Poincar\'e} and resp. {\it mixed Hodge polynomials}, as well as {\it $Rep_{\bC}(G)$-valued characters} as follows:
\begin{itemize}
\item {\it $Rep_{\bK}(G)$-valued Poincar\'e polynomials}:
$$P^G_{({c})}(X,\cM)(z):=\sum_{k} [H_{({c})}^{k}(X,\cM)]\cdot (-z)^k \in Rep_{\bK}(G)[z^{\pm 1}].$$
\item {\it $Rep_{\bC}(G)$-valued mixed Hodge polynomials}:
$$h^G_{({c})}(X,\cM)(y,x,z):=\sum_{p,q,k} [H_{({c})}^{p,q,k}(X,\cM)] \cdot y^px^q(-z)^k  \in Rep_{\bC}(G)[y^{\pm 1},x^{\pm 1},z^{\pm 1}].$$
\item {\it $Rep_{\bK}(G)$-valued equivariant characters}:
$$tr^G_{\Sigma_n}(H^*_{({c})}(X^n,\cM^{\boxtimes n})):=\sum_{k}tr^G_{\Sigma_n}(H_{({c})}^{k}(X^n,\cM^{\boxtimes n})) \cdot (-z)^k \in C(\Sigma_n) \otimes_{\bZ} Rep_{\bK}(G) \otimes_{\bZ} \bL,$$
with $\bL=\bZ[z^{\pm 1}]$, and, resp., $\bL=\bZ[y^{\pm 1}, x^{\pm 1}, z^{\pm 1}]$ in the Hodge context (with $\bK=\bC$ in this case).
\end{itemize}

 \br 
 If $G$ is the trivial (semi-)group, then $A_G={\rm Vect}_{\bK}$ is the category of finite-dimensional $\bK$-vector spaces, and $\dim: Rep_{\bK}(G) \overset{\simeq}{\lra} \bZ$, so the above $Rep_{\bK}(G)$-valued equivariant Poincar\'e and mixed Hodge polynomials, resp., $Rep_{\bK}(G)$-valued characters reduce in this case to the classical notions from Sections \ref{gs} and \ref{twist} of the Introduction.
 \er
 
Results analogous to those presented in Sections \ref{gs} and \ref{twist} can now be formulated for these modified notions of invariants in the $G$-equivariant context, with the corresponding Adams operations $$\psi_r:Rep_{\bK}(G) \otimes_{\bZ} \bL \to Rep_{\bK}(G) \otimes_{\bZ} \bL $$ defined as the tensor product of the Adams operations on the tensor factors (with $Rep_{\bK}(G)$ a pre-lambda ring as before).
Moreover, their proofs follow as before from  Theorems \ref{mti} and \ref{itwist} in the abstract context, but using the category $A_G$ in place of $A$.

Let us illustrate such formulae analogous to (\ref{mf1}) and (\ref{new2a}) in the $G$-equivariant context for the case of Poincar\'e polynomial invariants. Similar results for the mixed Hodge context, as well as various specializations of the variables are left to the reader. 
\bt\label{eqtl} Let $X$ be a space with a $G$-action as above, and let $\cM \in A_G(X)$ be a $G$-equivariant object in $A(X)$. Then:
\be\label{mf1l}
\sum_{n \geq 0} tr^G_{\Sigma_n}(H^*_{({c})}(X^n,\cM^{\boxtimes n})) \cdot t^n  
 =\exp \left( \sum_{r \geq 1}  p_r \otimes {\psi_r(P^G_{({c})}}(X,\cM)(z)) \cdot \frac{t^r}{r}\right) 
\ee
holds in the graded $\bQ$-algebra $Rep_{\bK}(G) \otimes_{\bZ}  \bQ[p_i, i \geq 1,z^{\pm 1}][[t]]$, and, respectively,
\be\label{new2al}
tr^G_{\Sigma_n}(H^*_{({c})}(X^n, V \otimes \cM^{\boxtimes n})) = \sum_{{\lambda=(k_1,k_2, \cdots) \vdash n} } \frac{p_{\lambda}}{z_{\lambda}} \chi_{\lambda}(V)  \otimes \prod_{r \geq 1} \left( \psi_r(P^G_{({c})}(X,\cM)(z)) \right)^{k_r},
\ee
holds in $Rep_{\bK}(G) \otimes_{\bZ}  \bQ[p_i, i \geq 1,z^{\pm 1}]$, for a given $V\in Rep_{\bQ}(\Sigma_n)$.
\et

In this equivariant context, new specializations of Theorem \ref{eqtl} are obtained 
by applying suitable  ring homomorphisms $sp: Rep_{\bK}(G) \to R$ to a commutative ring $R$. More concretely, in the three situations (A)-({C}) considered at the beginning of this section, examples of such specializations $sp: Rep_{\bK}(G) \to R$ are given as follows:
\begin{itemize}
\item[(A)] for a finite group $G$, we take the characters of $G$-representations, i.e., we apply the pre-lambda ring homomorphism $$tr_G:Rep_{\bK}(G) \lra C(G) \otimes_{\bZ} \bK, \ \ V \mapsto \left(g\mapsto trace_g(V)\right),$$ with Adams operations $\psi_r$ on $C(G) \otimes \bK$ given by $(\psi_r(\alpha))(g):=\alpha(g^r)$, for $\alpha \in C(G)$ and $g \in G$.
\item[(B)] for $G=\bZ$, with $g=1\in \bZ$ acting with finite order, assume for simplicity that $\bK$ is algebraically closed (or make a base change); then 
we have a pre-lambda ring isomorphism 
$$sp:Rep_{\bK}(G) \simeq \bZ[\widehat{\mu}],$$ with $\widehat{\mu}$ the abelian group of roots of unity in $\bK$ (with respect to multiplication), given by $[\chi_{\lambda}] \mapsto (\lambda)$, where $\chi_{\lambda}$ is the one-dimensional representation with $1 \in \bZ$ acting by multiplication with $\lambda$. The $r$-th Adams operations $\psi_r$  on $\bZ[\widehat{\mu}]$ is defined by $(\lambda) \mapsto (\lambda^r)$, for all $\lambda \in \widehat{\mu}$ (i.e., it is induced from the group homomorphism
$\widehat{\mu}\to \widehat{\mu}; \lambda \mapsto \lambda^r$ of the abelian group of roots of unity $(\widehat{\mu},\cdot)$).
\item[(C)] for the endomorphism category ${\rm End}_{\bK}$, we consider the usual  ring  homomorphism $$trace:{\bar K}_0({\rm End}_{\bK}) \lra \bK$$ defined by taking the trace of the endomorphism, with \be\label{trl} trace\left(\psi_r(g:V \to V)\right)=trace(g^r:V\to V).\ee
The identity (\ref{trl}) can be obtained as follows. First, by base change we can assume that $\bK$ is algebraically closed. Then we factor $trace$ through the projection from ${\bar K}_0({\rm End}_{\bK})$ to the usual Grothendieck group of the abelian tensor category ${\rm End}_{\bK}$, which is a pre-lambda ring homomorphism (cf. \cite{MS}[Lemma 2.1]). Finally, we reduce via short exact sequences to the case of one-dimensional representations (given by eigenspaces). Note that $trace$ is not a pre-lambda ring homomorphism.  Pre-lambda ring homomorphisms relevant to this situation are: the {\it characteristic polynomial}:  $$\lambda_t: {\bar K}_0({\rm End}_{\bK}) \lra W_{rat}(\bK):=\{P(t)/Q(t) \ | \ P(t),Q(t) \in 1+t\bK[t] \} \subset \bK(t),$$ $$[V,g] \mapsto \lambda_t(V,g):=\det(1+tg)=\sum_{i \geq 0} trace_{\Lambda^i g}(\Lambda^i V) \cdot t^i$$
given by the traces of the induced endomorphisms of the alternating powers of $V$,
and respectively, the {\it $L$-function}:
$$[V,g] \mapsto L(V,g)(t):=\det(1-tg)^{-1}=\sum_{i \geq 0} trace_{Sym^i g}(Sym^i V) \cdot t^i$$
given by the traces of the induced endomorphisms of the symmetric powers of $V$.
Here, $W_{rat}(\bK)$ is the subring of {\it rational} elements (as in \cite{Na}[Prop.6]) in the {\it big Witt ring} 
$W(\bK):=(1+t\bK[[t]],\cdot)$, 
with a suitable ring structure as in \cite{Al, Ha, Ra}, and whose underlying additive structure is the multiplication of rational functions, resp., normalized formal powers series.
\end{itemize}

For simplicity, we illustrate some of these specializations of Theorem \ref{eqtl} in the case of  constant (Hodge) coefficients in the complex quasi-projective context (with all actions under consideration being algebraic).
Due the algebraic nature of the action, the (compactly supported) cohomology $H^*_{({c})}(X,\bQ)$ gets an induced pullback  action of: (A) the finite group $G$;  (B)  
the cyclic group $\langle g \rangle$, and, resp., (C) the endomorphism $g$, compatible with the Deligne mixed Hodge structures (with the assumption that $g$ is proper if $H_c(-)$ is considered). It follows that the graded pieces $H^{p,q,k}_{({c})}(X,\bC)$ carry a similar action. So we can define a corresponding {\it equivariant mixed Hodge polynomial} $h^G_{({c})}(X)$, $h^{\langle g \rangle}_{({c})}(X)$, and resp.  $h^g_{({c})}(X)$ in this equivariant context as follows:
\begin{itemize}
\item[(A)] If $G$ is a finite group,
$$h^G_{({c})}(X)(y,x,z):=\sum_{p,q,k} tr_G(H_{({c})}^{p,q,k}(X,\bC)) \cdot y^px^q(-z)^k  \in C(G) \otimes_{\bZ} \bC[y,x,z],$$
with $C(G) \otimes_{\bZ}  \bC$ the complex valued class-functions of $G$, and $tr_G$ the usual character map.
\item[(B)] If $g$ is an algebraic automorphism of $X$ of finite order,
we let 
$$\chi_{\langle g \rangle}^{p,q,k}(X):=\sum_{\lambda \in \widehat{\mu}} \dim_{\bC}(H_{({c})}^{p,q,k}(X,\bC)_{\lambda}) \cdot (\lambda) \in \bZ[\widehat{\mu}],$$
with $\bZ[\widehat{\mu}]$ the group ring of the abelian group $\widehat{\mu}$ of roots of unity in $\bC$ (with respect to multiplication), and $H_{({c})}^{p,q,k}(X,\bC)_{\lambda}$ denoting the corresponding $\lambda$-eigenspace of $g$. Then we set
$$h^{\langle g \rangle}_{({c})}(X)(y,x,z):=\sum_{p,q,k} \chi_{\langle g \rangle}^{p,q,k}(X) \cdot y^px^q(-z)^k  \in  \bZ[\widehat{\mu}] \otimes_{\bZ}  \bC[y,x,z].$$
\item[(C)] If $g:X \to X$ is a (proper) algebraic endomorphism, then we set
$$h^g_{({c})}(X)(y,x,z):=\sum_{p,q,k} trace_g(H_{({c})}^{p,q,k}(X;\bC)) \cdot y^px^q(-z)^k  \in  \bC[y,x,z].$$
\end{itemize}

The external products $X^n$ get an induced diagonal action of $G$,  $\langle g \rangle$ or, resp., $g$, commuting with the symmetric group action. Therefore, the symmetric products $X^{(n)}$ inherit a similar action of $G$, $\langle g \rangle$  or, resp., $g$, so the corresponding invariants from above are also defined for each $X^{(n)}$. Therefore, the isomorphism
$$H^*_{({c})}(X^n,\bQ)^{\Sigma_n} \simeq H^*_{({c})}(X^{(n)},\bQ)$$
is $G$-equivariant.

By further specializing in Theorem \ref{eqtl} all $p_r$'s to the value $1$ (correponding to taking the $\Sigma_n$-invariant part), one obtains via the above ring homomorphisms 
 the following Macdonald-type generating series result:
 \bc\label{eqt} 
\begin{itemize}
\item[(A)] If $G$ is a finite group acting algebraically on $X$, then:
\be\label{eqta}
 \sum_{n \geq 0}  h^G_{({c})}(X^{(n)}) \cdot t^n  
=\exp \left( \sum_{r \geq 1}  \psi_r(h^G_{({c})}(X)) \cdot \frac{t^r}{r} \right) \in C(G) \otimes_{\bZ} \bC[y,x,z][[t]],
\ee
with $\psi_r(h^G_{({c})}(X)(y,x,z))(g):=h^G_{({c})}(X)(y^r,x^r,z^r)(g^r)$, for all $g \in G$.
\item[(B)] If $g$ is an algebraic automorphism of $X$ of finite order, then:
\be\label{eqtb}
 \sum_{n \geq 0}  h^{\langle g \rangle}_{({c})}(X^{(n)}) \cdot t^n  
=\exp \left( \sum_{r \geq 1}  \psi_r(h^{\langle g \rangle}_{({c})}(X)) \cdot \frac{t^r}{r} \right) \in \bZ[\widehat{\mu}] \otimes_{\bZ}  \bC[y,x,z][[t]],
\ee
with $\psi_r((\lambda) \cdot h(y,x,z)) := (\lambda^r) \cdot h(y^r,x^r,z^r)$, for $\lambda \in \widehat{\mu}$ and $h(y,x,z) \in \bC[y,x,z]$.
\item[(C)] If $g:X \to X$ is a (proper) algebraic endomorphism of $X$, then 
\be\label{eqtc}
 \sum_{n \geq 0}  h^g_{({c})}(X^{(n)})(y,x,z) \cdot t^n  
=\exp \left( \sum_{r \geq 1}  h^{g^r}_{({c})}(X)(y^r,x^r,z^r) \cdot \frac{t^r}{r} \right) \in  \bC[y,x,z][[t]].
\ee
\end{itemize}
\ec
Let us next compare special cases of Corollary \ref{eqt} with other results available in the literature. 
\begin{itemize}
\item[(A)] By specializing to $z=1$, our invariant $h^G_{({c})}$ becomes the corresponding {\it equivariant $E$- (or Hodge-Deligne) polynomial} $E^G_{({c})}$.  By further specializing also $y$ and $x$ to the value $1$, this reduces to 
the more classical {\it equivariant Euler characteristic} $\chi^G_{({c})} \in C(G) \otimes_{\bZ}  \bC$. Then (\ref{eqta}) becomes a variant of \cite{GLM}[Lemma 1], which is formulated in terms of the Burnside ring $A(G)$ of $G$, instead of class functions.
\item[(B)] By specializing to $z=1$, our invariant $h^{\langle g \rangle}_{({c})}$ becomes the corresponding {\it equivariant $E$- (or Hodge-Deligne) polynomial} $E^{\langle g \rangle}_{({c})}$. Then formula (\ref{eqtb}) reduces in case of compact supports to \cite{EG}[Theorem 1], which is formulated in terms of the power structure on the pre-lambda ring $\bZ[\widehat{\mu}] \otimes_{\bZ}  \bC[y,x]$. By further specializing to $x=1$, this equivariant $E$-polynomial reduces to the well-studied  {\it  Hodge spectrum} of a finite order automorphism.
\item[(C)] The right-hand side of formula (\ref{eqtc}) is a Hodge version of the graded  {\it Lefschetz Zeta function}, to which it reduces by specializing the variables $y$ and $x$ to the value $1$. More precisely, (\ref{eqtc}) specializes to formula (\ref{eqf}) of Theorem \ref{teq}, with $\bK=\bQ$. The classical Lefschetz Zeta function is obtained by further specializing to $z=1$. Similarly, as $g^r=id_X$ for all $r$ in case $g=id_X$,  the graded (resp., Hodge) version of the classical  {Lefschetz Zeta function} specializes in this case to  (Cheah's  Hodge version \cite{Che} of) 
{Macdonald's generating series formula} \cite{Mac} for the Poincar\'e polynomials and Betti numbers of the symmetric products of $X$ (see formula (\ref{Mc1final}) and also Corollary \ref{final} below for the corresponding graded version of the Lefschetz Zeta function).
\end{itemize}

\br
The interested reader should compare our results also with \cite{Ka}[Prop.15.5] and, resp., \cite{Bi2}[Thm.3.12], for an abstract analog of (\ref{eqtc}) and, resp., (\ref{eqta})
in the context of an automorphism resp., of a finite group action for a {\it dualizable} object in a suitable tensor category, with a corresponding notion of a trace.
\er

The specialization at $z=1$ of Corollary \ref{eqt} (A), resp., (B) above can also be reformulated (compare also with  \cite{EG, GLM}) by saying that
\be
E^G_{c}: K^G_0(var/\bC)\to C(G) \otimes_{\bZ}  \bC[y,x] \quad \text{resp.} \quad 
E^{\langle g \rangle}_{c}: K^{\langle g \rangle}_0(var/\bC)\to  \bZ[\widehat{\mu}] \otimes_{\bZ}  \bC[y,x]
\ee
is a morphism of pre-lambda rings, with the  pre-lambda structure of  the corresponding equivariant Grothendieck group of complex algebraic varieties (with respect to the scissor relation) defined via the {\it Kapranov Zeta function}
$$[X]\mapsto [pt]+\sum_{n\geq 1} [X^{(n)}]\cdot t^n\:.$$
Similar considerations apply for the variant
\be
h^G_{({c})}:\bar{K}^G_0(var/\bC)\to C(G) \otimes_{\bZ}  \bC[y,x,z] \quad \text{resp.} \quad 
h^{\langle g \rangle}_{({c})}: \bar{K}^{\langle g \rangle}_0(var/\bC)\to  \bZ[\widehat{\mu}] \otimes_{\bZ}  \bC[y,x,z]
\ee
on the corresponding equivariant Grothendieck group of complex algebraic varieties (with respect to disjoint unions), as studied in \cite{MS}[Sec.2.2] in the non-equivariant context.

\medskip

As a final example, let us formulate the {\it graded version of the classical Lefschetz Zeta function} (corresponding to the use of the trace homomorphism in the context (C) as above, for the constant constructible sheaf, and with all Fobenius parameters $p_r=1$), see Theorem \ref{teq} from the Introduction:

\bc\label{final} 
If $g:X \to X$ is a (proper) algebraic endomorphism of $X$, then the following equalities hold in $\bK[z][[t]]$:
\be\begin{split}\label{efinal}
 \sum_{n \geq 0}  P^g_{({c})}(X^{(n)})(z) \cdot t^n  
&=\exp \left( \sum_{r \geq 1}  P^{g^r}_{({c})}(X)(z^r) \cdot \frac{t^r}{r} \right) \\
&=\exp \left(\sum_{k\geq 0} (-1)^k\left( \sum_{r \geq 1}  trace_{g^r}(H^k_{({c})}(X,\bK)) \cdot \frac{(z^kt)^r}{r} \right)\right) \\
&= \prod_{k\geq 0}\left(  \sum_{i \geq 0} trace_{Sym^i g}(Sym^i (H^k_{({c})}(X,\bK))) \cdot (z^kt)^i  \right)^{(-1)^{k}}\\
&= \prod_{k\geq 0}\left( L(H^k_{({c})}(X,\bK),g)(z^kt)\right)^{(-1)^{k}}
\end{split}\ee
\ec

Note that formula (\ref{efinal}) specializes for $z=1$ to the classical {Lefschetz Zeta function} of the (proper) endomorphism
$g: X\to X$:
\be\label{Leffinal}
\sum_{n \geq 0}  \chi^g_{({c})}(X^{(n)})(z) \cdot t^n  = \prod_{k\geq 0}\left( L(H^k_{({c})}(X,\bK),g)(t)\right)^{(-1)^{k}}\:.
\ee
On the other hand, for $g=id_X$ the identity of $X$, formula (\ref{efinal})
reduces to {Macdonald's generating series formula} \cite{Mac} for the Poincar\'e polynomials and Betti numbers of the symmetric products of $X$:
\be\label{Mc1final}
 \sum_{n \geq 0}  P_{({c})}(X^{(n)})(z) \cdot t^n =\exp \left( \sum_{r \geq 1}  P_{({c})}(X)(z^r) \cdot \frac{t^r}{r} \right)
= \prod_{k\geq 0}\left( \frac{1}{1-z^kt}\right)^{(-1)^{k}\cdot b^k_{({c})}(X)} \:,
\ee
with $b^k_{({c})}(X):=\dim_{\bC}H^k_{({c})}(X,\bK)$,
which for $z=1$ specializes (also as the particular case of (\ref{Leffinal}) for $g=id_X$) to:
\be\label{Mc2final}
 \sum_{n \geq 0}  \chi_{({c})}(X^{(n)})\cdot t^n 
= {(1-t)}^{-\chi_{({c})}(X)} \:.
\ee

\br 
For the counterpart of (\ref{Leffinal})
in the context of the
Zeta function of a constructible sheaf for the Frobenius endomorphism of varieties over finite fields, see also  \cite{SGA5}[Thm. on p.464] and \cite{FK}[Thm.4.4 on p.174]. 
For a similar counterpart of (\ref{efinal}) taking  a weight filtration into account, see \cite{Na}[Prop.8(i)].

Finally, the product on the big Witt ring $W(\bK)$ (or its subring $W_{rat}(\bK)$) corresponds under the ring homomorphisms
$\lambda_t, L(t): {\rm End}_{\bK}\to W_{rat}(\bK)\subset W(\bK)$ to the tensor product of endomorphisms. By the specialization above  of the (graded version of the) Lefschetz Zeta function to Macdonald's generating series formula for the Poincar\'e polynomials and Betti numbers of the symmetric products of $X$, it should not come as a surprise that the Witt multiplication naturally arises if one attempts to express these generating functions for a product space $X\times X'$ in terms of the corresponding 
generating functions of the factors $X$ and $X'$ (as further discussed in \cite{Ra, RT}).
\er

%%%%%%%%%%%%%%%%%%%%%%%%%%%%%%
%%%%%%%%%%%%%%%%%%%%%%%%%%%%%%

\section{Pseudo-functors and applications}\label{psf}

In this final section, we explain the connection of Theorem \ref{mth} with our previous results from \cite{MS} about generating series of symmetric and alternating powers of suitable coefficients, e.g., (complexes of) constructible or coherent sheaves, or (complexes of) mixed Hodge modules in the mixed Hodge context. 
In fact, all of this can be discussed in the abstract setting of suitable pseudo-functors, as in \cite{MS}, which we now recall. For simplicity, we focus here on the complex quasi-projective context. Note that this pseudo-functoriality of our coefficients as above is also needed if one works, as in \cite{MS15}, with the corresponding results for equivariant characteristic classes of these coefficients.

\subsection{Pseudo-functors}\label{psfs}

Let $(-)_*$ be a (covariant) pseudo-functor on the category of complex quasi-projective varieties (with proper morphisms), taking values in a
pseudo-abelian (or Karoubian) $\bQ$-linear additive category $A(-)$, e.g., see \cite{MS}[Sect.4.1]. In fact, our abstract axiomatic approach would also work for a suitable (small) category of  {\it spaces} with finite products and a terminal object $pt$ (corresponding to the empty product, see \cite{MS}[Appendix] for more details).
 Assume, moreover, that the following properties are satisfied:
\begin{enumerate}
 \item[(i)] For any quasi-projective variety $X$ and all $n$ there is a multiple external product 
$$\boxtimes^n: \;A(X)^{ \times n} \to A(X^n),$$ equivariant with respect to a permutation action of the symmetric group $\Sigma_n$,
i.e., $M^{\boxtimes n}\in A(X^n)$ is a $\Sigma_n$-equivariant object, for all $M\in A(X)$.
\item[(ii)] $A(pt)$ is endowed with a $\bQ$-linear tensor structure $\otimes$, which makes it into a symmetric monoidal category.
\item[(iii)]  For any quasi-projective variety $X$, $M\in A(X)$ and all $n$, there is a $\Sigma_n$-equivariant isomorphism
$$k_*(M^{\boxtimes n})\simeq (k_*M)^{\otimes n},$$ with $k$  the constant morphism to a point $pt$. Here, the $\Sigma_n$-action on the left-hand side
is induced from (i), whereas the one on the right-hand side comes from (ii).
\end{enumerate}
For example, the above properties are fullfilled for $A(X)=D^b\mh(X)$, the bounded derived category of algebraic mixed Hodge modules on $X$, viewed as a pseudo-functor with respect to either of the push-forwards $(-)_*$ or $(-)_!$, 
as well as for the derived categories $D^b_c(X;\bK)$ and $D^b_{coh}(X)$ of bounded complexes with constructible and, resp., coherent cohomology, see \cite{MS} for more details. In the coherent setting, we restrict to projective varieties $X$, so that in this context $(-)_*=(-)_!$. As another example, one can take $A(X)=D^b_{f}(X;\bK)$ to be the triangulated category consisting of bounded complexes of sheaves of $\bK$-vector spaces with finite dimensional compactly supported cohomology, viewed as a pseudo-functor with respect to $(-)_!$.

\br
Property (iii) is the abstract analogue of the K\"unneth isomorphism (\ref{Kue}).
\er

Let $\pi_n : X^n\to X^{(n)}$ be the natural projection onto the $n$-th symmetric product $X^{(n)}:=X^n/\Sigma_n$. By property (i), for any $M \in A(X)$ the exterior product $M^{\boxtimes n}$ is a $\Sigma_n$-equivariant object in $A(X^n)$, i.e., it is an element of $A_{\Sigma_n}(X^n)$, e.g., see \cite{MS}[Sect.4.2]. Then the pushdown ${\pi_n}_*M^{\boxtimes n}$ to the $n$-th symmetric product is a $\Sigma_n$-equivariant object on $X^{(n)}$. Since $\Sigma_n$ acts trivially on $X^{(n)}$, the $\Sigma_n$-action on ${\pi_n}_*M^{\boxtimes n}$ corresponds to a group homomorphism $$\Psi:\Sigma_n \to \Aut_{A(X^{(n)})}({\pi_n}_*M^{\boxtimes n}).$$ Moreover, since $A(X^{(n)})$ is a $\bQ$-linear additive category, we can define the symmetric projector $$(-)^{\Sigma_n}:=\frac{1}{n!} \sum_{\sigma \in \Sigma_n} \Psi_{\sigma}$$ onto the $\Sigma_n$-invariant part, and, respectively, the alternating projector $$(-)^{sign-\Sigma_n}:=\frac{1}{n!} \sum_{\sigma \in \Sigma_n} {sign(\sigma)} \Psi_{\sigma}, $$ for $sign: \Sigma_n \to \{\pm 1\}$ the sign character, and $\Psi_{\sigma}$ denoting the $\sigma$-action $\Psi(\sigma)$.
Using the Karoubian structure, we can then associate to an object $M \in A(X)$ its {\it $n$-th symmetric power} $$M^{(n)}:=\left( \pi_{n*}M^{\boxtimes n} \right)^{\Sigma_n}$$
and, respectively, its {\it $n$-th alternating power} $$M^{\{n\}}:=\left( \pi_{n*}M^{\boxtimes n} \right)^{sign-\Sigma_n},$$ as objects in $A(X^{(n)})$.  As in \cite{MS}[Sect.2], we then have the identities (with $k$ denoting in this paper the constant map from any space to a point):
\begin{equation}\label{key}
k_*(M^{(n)}) \simeq \left((k_*M)^{\otimes n}\right)^{\Sigma_n}\simeq (k_*M)^{(n)} \quad \text{and}  \quad
k_*(M^{\{n\}}) \simeq \left((k_*M)^{\otimes n}\right)^{sign-\Sigma_n}\simeq (k_*M)^{\{n\}} ,
\end{equation}
which allow the calculation of invariants of $k_*M^{(n)}$ and $k_*M^{\{n\}}$, respectively, 
only in terms of those for $k_*M \in A(pt)$ and the symmetric monoidal structure $\otimes$, see \cite{MS} for more details.
Here we are interested in representation-theoretic refinements of such formulae from \cite{MS} expressed in terms of abstract generating series identities  for the $\Sigma_n$-equivariant objects ($n \geq 0$): 
$$k_*M^{\boxtimes n}\simeq (k_*M)^{\otimes n} \in A_{\Sigma_n}(pt)\:.$$

\medskip

In this section $A(pt)=:A$ plays the role of the underlying symmetric monoidal category used in Section \ref{simmon}.

\medskip

Let $\bar{K}_0(-)$ denote as before the Grothendieck group of an additive
category viewed as an exact category by the split exact sequences corresponding to direct sums
$\oplus$, i.e., the Grothendieck group associated to the abelian monoid of isomorphism classes of objects with the direct sum. 
As in Section \ref{simmon}, $\bar{K}_0(A(pt))$ becomes a pre-lambda ring.

By Theorem \ref{mta}, applied to $$\cV:=k_*\cM \in A(pt),$$ with $M \in A(X)$, we obtain by property (iii) of the pseudo-functor $(-)_*$ the following  generalization of \cite{MS}[Thm.1.7]:
\bt\label{mtl} For any $M \in A(X)$, the following generating series identity holds in the $\bQ$-algebra $\left( \bar{K}_0(A(pt)) \otimes_{\bZ} \bQ[p_i, i \geq 1] \right)[[t]]  = \left( \bQ[p_i, i \geq 1]\otimes_{\bZ} \bar{K}_0(A(pt))  \right)[[t]] $: 
\be\label{ml}  \sum_{n \geq 0}  cl_n([k_*M^{\boxtimes n}]) \cdot t^n=\exp \left( \sum_{r \geq 1} p_r \otimes \psi_r([k_*M]) \cdot \frac{t^r}{r}    \right),
\ee
with $\psi_r$ the corresponding $r$-th Adams operation of the pre-lambda ring $\bar{K}_0(A(pt))$.
\et

Specializing to $p_r=1$ for all $r$ corresponds to the functor induced on Grothendieck groups by taking the $\Sigma_n$-invariant part 
$$(-)^{\Sigma_n}=\frac{1}{n!} \sum_{\sigma \in \Sigma_n} \Psi_{\sigma}:A_{\Sigma_n}(pt) \lra A(pt).$$ Indeed, this reduces via the decomposition
$$\bar{K}^{\Sigma_n}_0(A(pt)) \simeq \bar{K}_0(A(pt)) \otimes_{\bZ} Rep_{\bQ}(\Sigma_n)$$
to the corresponding classical formula for the represention ring $Rep_{\bQ}(\Sigma_n)$.
So, by letting $p_r=1$ for all $r$ in Theorem \ref{mtl}, one obtains by the isomorphism $$k_*(M^{(n)}) \simeq \left((k_*M)^{\otimes n}\right)^{\Sigma_n}\simeq (k_*M)^{(n)}$$ the following generating series from \cite{MS}[Thm.1.7]:
\be\label{sim2} 
1+\sum_{n\geq 1}\;  [k_*M^{(n)}] \cdot t^n =\exp\left( \sum_{r\geq 1}\;
\psi_r([k_*M]) \cdot \frac{t^r}{r}\right) \in \bar{K}_0(A(pt))\otimes_{\bZ}\bQ[[t]]  \:.\ee

\medskip

Similarly, specializing to $p_r=(-1)^{r-1}={sign(\sigma_r)}$ for all $r$ (with $\sigma_r$ denoting as before an $r$-cycle in $\Sigma_r$) corresponds to the functor induced on Grothendieck groups by taking the projector onto the alternating part of the $\Sigma_n$-action:
$$(-)^{sign-\Sigma_n}=\frac{1}{n!} \sum_{\sigma \in \Sigma_n} {sign(\sigma)} \Psi_{\sigma}:A_{\Sigma_n}(pt) \lra A(pt).$$ 
So, by letting $p_r=(-1)^{r-1}$ for all $r$ in Theorem \ref{mtl}, one obtains by the isomorphism $$k_*(M^{\{n\}}) \simeq \left((k_*M)^{\otimes n}\right)^{sign-\Sigma_n} \simeq (k_*M)^{\{n\}}$$ the following generating series from \cite{MS}[Thm.1.7]:
\be\label{alt2}1+\sum_{n\geq 1}\; [k_*M^{\{n\}}] \cdot t^n =\exp\left(- \sum_{r\geq 1}\;
\psi_r([k_*M]) \cdot \frac{(-t)^r}{r}\right) \in \bar{K}_0(A(pt))\otimes_{\bZ}\bQ[[t]]  \:.\ee

\medskip

We next explain how to derive our concrete formulae of Theorem \ref{mth} from the Introduction also by using Theorem \ref{mtl} of this section, by applying suitable pre-lambda ring homomorphisms. The virtue of this approach is that it explains the connection of Theorem \ref{mth} with our previous results from \cite{MS} about generating series of symmetric and alternating powers of suitable coefficients, and it also yields similar generating series identities on various intermediate levels of Grothendieck groups between $\bar{K}_0(D^b\mh(pt))$ and $\bZ[y^{\pm 1},x^{\pm 1},z^{\pm 1}]$.

Consider the homomorphism of pre-lambda rings
$$h:  \bar{K}_0(D^b\mh(pt)) \to \bZ[y^{\pm 1},x^{\pm 1},z^{\pm 1}]$$
defined via the commutative diagram of homomorphisms of pre-lambda rings (as in \cite{MS}[p.301]):

\begin{equation*}
\begin{CD}
 \bar{K}_0(D^b\mh(pt)) @> H^* >> \bar{K}_0(Gr^{-}(\mh(pt))) @> \sim >> \bar{K}_0(Gr^{-}(\ms^p)) \\
@V h VV   @VVV @VV {\rm forget} V \\
 \bZ[y^{\pm 1},x^{\pm 1},z^{\pm 1}] @<< h < \bar{K}_0(Gr^{-}(Gr^2({\rm Vect}_{\bC}))) @< Gr_F^*Gr^W_* << \bar{K}_0(Gr^{-}(\ms)) \:.  
  \end{CD}
\end{equation*}
The bottom row was already explained in Section \ref{pmt}. Additionally,
the following notations are used:
\begin{enumerate}
\item[(a)] $H^*: D^b\mh(pt)\to Gr^{-}(\mh(pt))$ is the total cohomology functor $\cV\mapsto \bigoplus_n H^n(\cV)$. Note that this is a functor of additive tensor categories
(i.e., it commutes with direct sums $\oplus$ and tensor products $\otimes$), if we choose the Koszul symmetry isomorphism on $Gr^{-}(\mh(pt))$.
In fact, $D^b\mh(pt)$ is a triangulated category with bounded $t$-structure satisfying \cite{Bi}[Def.4.2], so that the claim follows from \cite{Bi}[Thm.4.1, Cor.4.4].
\item[(b)] The isomorphism $\mh(pt)\simeq \ms^p$ is Saito's identification of the abelian tensor category of mixed Hodge modules over a point space with Deligne's abelian tensor category of polarizable mixed Hodge structures.
\item[({c})] ${\rm forget}: \ms^p\to \ms$ is the functor of forgetting that the corresponding $\bQ$-mixed Hodge structure is graded polarizable.
\end{enumerate}

\br\label{rem-kue}
The fact that the total cohomology functor $H^*: D^b\mh(pt)\to Gr^{-}(\mh(pt))$ is a tensor functor corresponds to the  K\"{u}nneth formula 
$$H^*(\cV^{\otimes n}) \simeq (H^*(\cV))^{\otimes n} \ , \quad \text{ for $\cV\in D^b\mh(pt)$.}$$
For $\cV=k_*M$, this
implies by Property (iii) the important K\"{u}nneth isomorphism (\ref{Kue}) from the Introduction. For a more direct approach to K\"{u}nneth formulae, see \cite{Sch}[eq.(1.17), Cor.2.0.4] and \cite{MSS}[Sect.3.8] for the constructible context, \cite{B}[Thm.2.1.2] for the coherent context and, resp., \cite{MSS}[Thm.1] for the mixed Hodge module context.
\er

%\medskip

Formula (\ref{mf1h}) follows now by applying the ring homomorphism 
$$h  \otimes id : {\bar K}_0(D^b\mh(pt)) \otimes_{\bZ} \bQ[p_i, i \geq 1] \lra \bZ[y^{\pm 1}, x^{\pm 1},z^{\pm 1}] \otimes_{\bZ} \bQ[p_i, i \geq 1]$$
to formula (\ref{ml}) of Theorem \ref{mtl}.

\medskip

By exactly the same method one also gets the following homomorphism of pre-lambda rings:
$$\begin{CD} 
 P: \bar{K}_0(D^b_c(pt)) @> H^* >>  \bar{K}_0(Gr^{-}({\rm Vect}_{\bK}))) @> P >> \bZ[z^{\pm 1}]
  \end{CD}$$
and, resp., 
$$\begin{CD} 
 P: \bar{K}_0(D^b_{coh}(pt)) @> H^* >>  \bar{K}_0(Gr^{-}({\rm Vect}_{\bK}))) @> P >> \bZ[z^{\pm 1}]\:,
  \end{CD}$$
with $P: \bar{K}_0(Gr^{-}({\rm Vect}_{\bK})))\to \bZ[z^{\pm 1}]$ the Poincar\'e polynomial homomorphism defined in Section \ref{pmt}.
Then formula (\ref{mf1}) follows by applying  
$$P  \otimes id : {\bar K}_0(A(pt)) \otimes_{\bZ} \bQ[p_i, i \geq 1] \to \bZ[z^{\pm 1}] \otimes_{\bZ} \bQ[p_i, i \geq 1]$$
to formula (\ref{ml}), 
where $A(pt)$ is either $D^b_c(pt;\bK)=D^b_f(pt;\bK)$ or $D^b_{coh}(pt)$.

%%%%%%%%%%%%%%%%%%%%%%%%%%%%%%%%%%
%%%%%%%%%%%%%%%%%%%%%%%%%%%%%%%%%%

\subsection{Pseudo-functors and twisting}\label{psftwist} In the context of twisting by representations, we need to require the pseudo-functor $(-)_*$ with values in the category $A(-)$ to satisfy an additional property:
\begin{itemize}
\item[(iv)] For any quasi-projective variety $X$, there exists a pairing 
$$\otimes: A(pt) \times A(X) \lra A(X),$$
which is additive, $\bQ$-linear and functorial in each variable, as well as functorial with respect to $(-)_*$. Moreover, if $X=pt$ is a point, this pairing coincides with the tensor structure on $A(pt)$ of property (ii).
\end{itemize}
The pairing of (iv) induces similar ones on the corresponding equivariant categories, as well as on the (equivariant) Grothendieck groups. These pairings are bilinear and functorial with respect to the pseudo-functor $(-)_*$. Note that this additional property is fullfilled for all examples of pseudo-functors considered in this paper, i.e., $A(X)=D^b\mh(X)$, $D^b_c(X;\bK)$, $D^b_{coh}(X)$, or $D^b_f(X;\bK)$, where it is given as a special case of the exterior product $\boxtimes$, with 
$$\otimes : = k_*(-\boxtimes -): A(pt)\times A(pt)\to A(pt) $$
for $k: pt\times pt\simeq pt$. As before, in the coherent setting we restrict to projective varieties $X$. 

\br\label{ptwist}
In the context of our examples, the category ${\rm Vect}_{\bQ}(\Sigma_n)$ is a tensor subcategory of $A_{\Sigma_n}(pt)$, where in the Hodge context we regard a representation as a pure Hodge structure of type $(0,0)$ placed in degree zero, together with Saito's identification $\mh(pt)\simeq \ms^p$. Property (iv) yields now the pairing mentioned in Sect.\ref{twist}: $$\otimes: {\rm Vect}_{\bQ}(\Sigma_n) \times A_{\Sigma_n}(X) \to A_{\Sigma_n}(X),$$
which is induced from the composition:
$$A_{\Sigma_n}(pt) \times  A_{\Sigma_n}(X) \overset{\otimes}{\lra} A_{\Sigma_n \times \Sigma_n}(X) \overset{{\rm Res}}{\lra} A_{\Sigma_n}(X),$$
with ${\rm Res}$ the restriction functor for the diagonal subgroup $\Sigma_n \subset \Sigma_n \times \Sigma_n$. Moreover, if $X=pt$ is a point space, this pairing coincides with the abstract pairing (\ref{apair}) defined via the $\bQ$-linear structure.
\er

\br By the functoriality of the above pairing, we have the following  projection formula for a morphism $f:X \to X'$, $V \in {\rm Vect}_{\bQ}(\Sigma_n)$ and $\cM \in A_{\Sigma_n}(X)$:
\be\label{projf} f_*(V \otimes \cM)=V \otimes f_*(\cM),
\ee
where we also use  the identification $id_{pt} \times f=f$.
Applying this formula for $f$ the constant map $k:X^n \to pt$, together with the tensor property of the total cohomology functor $H^*$ as in Remark \ref{rem-kue}, we get the first isomorphism of the equivariant K\"unneth formula (\ref{tKue}).
\er

\bd For $V \in {\rm Vect}_{\bQ}(\Sigma_n)$ a rational $\Sigma_n$-representation, 
the {\it Schur-type object} $S_{V}(\cM) \in A(X^{(n)})$ associated to $\cM \in A(X)$ is defined by
\be
S_{V}(\cM):=\left( V \otimes \pi_{n*}(M^{\boxtimes n}) \right)^{\Sigma_n}.
\ee
If $V=V_{\mu} \simeq V^*_{\mu}$ is the (self-dual) irreducible representation of $\Sigma_n$ corresponding to a partition $\mu$ of $n$, we denote the corresponding Schur functor by $S_{\mu}:=S_{V_\mu}$.
\ed

Note that for $V$ the trivial (resp., sign) representation of $\Sigma_n$, the corresponding Schur functor coincides with the symmetric (resp., alternating) $n$-th power of $\cM$. Moreover, by using the projection formula for the constant map $k$ to a point, we have that
$$k_* S_V(\cM):=k_*\left( V \otimes \pi_{n*}(\cM^{\boxtimes n}) \right)^{\Sigma_n} \simeq 
\left( V \otimes k_*(\cM^{\boxtimes n}) \right)^{\Sigma_n}\simeq S_V(k_*\cM),$$
with the last identification following from Property (iii) of the pseudo-functor $(-)_*$.
Together with the tensor property of the total cohomology functor $H^*$ as in Remark \ref{rem-kue}, this yields formula (\ref{sc}). 

\bex
Let us now revisit the Schur functor $S_V$ obtained by choosing $V={\rm Ind}_K^{\Sigma_n}(triv_K)$, the representation induced from the trivial representation of a subgroup $K$ of $\Sigma_n$.  Then, if $\pi:X^n \lra X^n/K$ and $\pi':X^n/K \to X^{(n)}$ are the projections factoring $\pi_n$, we have:
$$ (\pi_{n*} (V \otimes \cM^{\boxtimes n}))^{\Sigma_n}
\simeq ({\rm Ind}_K^{\Sigma_n}(triv_K) \otimes \pi_{n*} (\cM^{\boxtimes n}))^{\Sigma_n}
\simeq (\pi_{n*} (\cM^{\boxtimes n}))^{K} \simeq \pi'_*\left(  (\pi_{*} (\cM^{\boxtimes n}))^{K}  \right),$$
for $\cM \in A(X)$.
As an example, if $\cM=\bQ^H_X \in D^b\mh(X)$ is the constant Hodge module on $X$, we get 
\be\label{hquot}
h_{({c})}(X^n/K)(y,x,z)= \sum_{{\lambda=(k_1,k_2, \cdots) \vdash n} } \frac{1}{z_{\lambda}} \chi_{\lambda}({\rm Ind}_K^{\Sigma_n}(triv))  \cdot \prod_{r \geq 1} \left( h_{({c})}(X)(y^r,x^r,z^r) \right)^{k_r},
\ee
and similarly for the Poincar\'e polynomials as in (\ref{pquot}).
\eex

%%%%%%%%%%%%%%%%%%%%%%%%%%%%%%
%%%%%%%%%%%%%%%%%%%%%%%%%%%%%%

\subsection{Pseudo-functors in the equivariant context}
Let $X$ be a complex quasiprojective variety, and $G$ be as in the beginning of Section \ref{feq}, i.e., (A) a finite algebraic group $G$, (B) a finite order algebraic automorphism $g:X \to X$, or (C) a (proper) algebraic endomorphism $g:X \to X$. For a $G$-equivariant object  $\cM \in A(X)$ as above, the external products $\cM^{\boxtimes n} \in A(X^n)$, their pushforwards $\pi_{n*}(\cM^{\boxtimes n})\in A(X^{(n)})$, and the symmetric and alternating powers $\cM^{(n)}$, $\cM^{\{n\}}  \in A(X^n)$ 
get an induced diagonal $G$-action commuting with the action of the symmetric group $\Sigma_n$ as before, so that for $V$ a $\Sigma_n$-representation (with trivial $G$-action), $V \otimes \cM^{\boxtimes n}$,  $\pi_{n*}(V \otimes \cM^{\boxtimes n})$ and  the Schur objects $S_V(\cM)$ get an induced $G$-action.

Therefore, all results from Sections \ref{psfs} and \ref{psftwist} hold in this equivariant context, provided that the derived K\"unneth formula of Property (iii) holds $G$-equivariantly.
In the abstract context of a pseudo-functor, this $G$-equivariance of the derived K\"unneth formula can be formulated as the following property of the pseudo-functor $(-)_*$ with values in the category $A(-)$:
\begin{itemize}
\item[(v)] For $g:X \to X$ an algebraic (iso)morphism and $\cM \in A(X)$ with a(n) (iso)morphism $\Psi_g: \cM \lra g_*\cM$ given by the $G$-action, we have an isomorphism 
\be
(g^{\times n})_*(\cM^{\boxtimes n}) \simeq (g_*\cM)^{\boxtimes n}
\ee
such that the (iso)morphism 
$$k_*\Psi_g^{\boxtimes n}: k_*\cM^{\boxtimes n} \to k_*\cM^{\boxtimes n}$$
induced by pushing down to a point (via $k_*$) the (iso)morphism
$$\Psi_g^{\boxtimes n}: \cM^{\boxtimes n} \lra (g_*\cM)^{\boxtimes n} \simeq (g^{\times n})_*(\cM^{\boxtimes n})$$
agrees under the identification $k_*\cM^{\boxtimes n} \simeq (k_*\cM)^{\otimes n}$ of Property (iii) with the endomorphism $$(k_*\Psi_g)^{\otimes n}: (k_*\cM)^{\otimes n} \lra (k_*\cM)^{\otimes n}.$$
In the case (A) of a finite group action, we ask this compatibility for all $g\in G$ (in such a way  that the corresponding $G$-actions via
$k_*\Psi_g^{\boxtimes n}$ and $(k_*\Psi_g)^{\otimes n}$ are identified under (iii)).
\end{itemize}
This Property (v) holds in the four main situations of coefficients considered here (i.e., $A(X)=D^b\mh(X)$, $D^b_c(X;\bK)$, $D^b_{coh}(X)$, or $D^b_f(X;\bK)$) by the {\it equivariance} of the corresponding  multiple K\"{u}nneth formula (as in the following Remark),
 see \cite{MS}[Appendix] for the constructible, coherent, and finite-dimensional context, and \cite{MSS}[Sect.1.12] for the Hodge context. 

\br The required compability (v) follows from the {\it equivariance} resp. {\it functoriality} of the multiple K\"{u}nneth formula:
$$(k^{\times n})_*(\boxtimes_{i=1}^n (-)) = \boxtimes_{i=1}^n (k_*(-)) : A(X)^{\times n} \to A(pt^{\times n})\:,$$
together with
$$\otimes_{i=1}^n (-) =   k_*(\boxtimes_{i=1}^n (-))  : A(pt)^{\times n} \to A(pt)\:.$$
The corresponding $G$-equivariance in the twisting defined via Property (iv) follows already from the required functorialities. 
\er

%%%%%%%%%%%%%%%%%%%%%%%%%%%%%%

\end{document}